%% file: SocialDynamicsSurvey.tex
\documentclass[graybox]{svmult}

\usepackage{mathptmx}       % selects Times Roman as basic font
\usepackage{helvet}         % selects Helvetica as sans-serif font
\usepackage{courier}        % selects Courier as typewriter font
\usepackage{type1cm}        % activate if the above 3 fonts are
\usepackage{makeidx}         % allows index generation
\usepackage{graphicx}        % standard LaTeX graphics tool
\usepackage{multicol}        % used for the two-column index
\usepackage[bottom]{footmisc}% places footnotes at page bottom

\usepackage{enumerate,stmaryrd,mathtools,amssymb}
\usepackage{float}
\usepackage{caption}
\usepackage{subfig}
\usepackage{tikz}
\usepackage{verbatim}
\usepackage{enumitem}
\newtheorem{prop}{Proposition}
\newtheorem{conj}{Conjecture}

\newcommand{\R}{\mathbb{R}}
\newcommand{\N}{\mathbb{N}}
\newcommand{\Sp}{\mathbb{S}}
\newcommand{\T}{\mathbb{T}}
\newcommand{\NN}{\mathcal{N}}

\newcommand{\GG}{\mathcal{G}}
\newcommand{\VV}{\mathcal{V}}
\newcommand{\EE}{\mathcal{E}}
\DeclareMathOperator{\tr}{tr}
\DeclareMathOperator{\dv}{div}

\makeindex             % used for the subject index
\begin{document} 

\date{\today}
\title{\LARGE \bf Interaction Network, State Space and Control in Social Dynamics}
\titlerunning{Interaction Network, State Space and Control}

\author{Aylin Aydo\u{g}du, Marco Caponigro, Sean McQuade, Benedetto Piccoli, Nastassia Pouradier Duteil, Francesco Rossi and Emmanuel Tr\'elat}
\authorrunning{Aydo\u{g}du, Caponigro, McQuade, Piccoli, Pouradier Duteil, Rossi, Tr\'elat }
%\author[*]{Aylin Aydogdu}
%\author[$\dag$]{Marco Caponigro}
%\author[*]{Maria Laura Delle Monache}
%\author[*]{Sean McQuade}
%\author[*]{Benedetto Piccoli}
%\author[*]{Nastassia Pouradier Duteil}
%\author[$\dag$]{Francesco Rossi}
%\author[$\dag$]{Emmanuel Tr\'elat}
%\affil[*]{ Department of Mathematical Sciences, Rutgers University, Camden, NJ 08102, USA.}%
%\affil[$\dag$]{Other affiliations}%

\institute{Aylin Aydo\u{g}du \at Rutgers University - Camden, Camden, NJ, \email{aylinvet87@gmail.com}
\and Marco Caponigro \at Conservatoire National des Arts et M\'etiers, Equipe M2N, 292 rue Saint-Martin, 75003, Paris, France, \email{marco.caponigro@cnam.fr}
\and Sean McQuade \at Rutgers University - Camden, Camden, NJ, \email{sean.mcquade@rutgers.edu}
\and Benedetto Piccoli \at Rutgers University - Camden, Camden, NJ, \email{piccoli@camden.rutgers.edu }
\and Nastassia Pouradier Duteil \at Rutgers University - Camden, Camden, NJ, \email{nastassia.pouradierduteil@rutgers.edu}
\and Francesco Rossi \at Aix Marseille Universit\'e, CNRS, ENSAM, Universit\'e de Toulon, LSIS UMR 7296,13397, Marseille, France, \email{francesco.rossi@lsis.org}
\and Emmanuel Tr\'elat \at Sorbonne Universit\'es, UPMC Univ Paris 06, CNRS UMR 7598, Laboratoire Jacques-Louis Lions, Institut
Universitaire de France, F-75005, Paris, France, \email{emmanuel.trelat@upmc.fr}}

\maketitle

%\tableofcontents

%%%%%%%%%%%%%%%%%%%%%%%%%%%%

\abstract{In the present chapter we study the emergence of global patterns in large groups in first and second-order multi-agent systems, focusing on two ingredients that influence the dynamics: the interaction network and the state space.
The state space determines the types of equilibrium that can be reached by the system. 
 Meanwhile, convergence to specific equilibria depends on the connectivity of the interaction network and on the interaction potential. 
When the system does not satisfy the necessary conditions for convergence to the desired equilibrium, control can be exerted, both on finite-dimensional systems and on their mean-field limit.
}

%%%%%%%%%%%%%%%%%%%%%%%%%%%%
\section{Introduction}
\input{Intro}

\section{Overview of Social Dynamics Problems}\label{Sec:Overview} % Sean, Aylin, Nastassia
\input{Overview}

\section{Role of the Interaction Network}
In finite-dimensional systems, the set of interacting agents can be interpreted as the vertices $\VV$ of a graph $\GG$, and their interactions can be represented as weighted edges $\EE$, as seen in Section \ref{Sec:Overview-gen}.  
In all the models reviewed in Section \ref{Sec:Overview}, the dynamics depend on the interaction network via the interaction coefficients $a_{ij}$ (see systems \eqref{eq:syst-gen} and \eqref{eq:syst2-gen}). 
In turn, the interaction network may depend on the dynamics, for instance when the interaction coefficients depend on the state variables: $a_{ij}=a(\|x_i-x_j\|)$. Then the graph $\GG$ varies in time. 
In this section we explore the influence of the network on the dynamics, and vice-versa. 

We will look at models in which (at least initially) the set of neighbors for each agent is smaller than the set of all agents: $\mathbf{card}(\NN_i)<N$, so $\EE\subsetneq \VV\times\VV $, such as bounded confidence models, as defined in Section \ref{Sec:Overview-order1}.

On the other hand, some models use the complete set of agents as the interaction network, so that each agent interacts with all the others. The interaction network is then interpreted as a weighted graph, where each edge's weight is given by the interaction coefficient $a_{ij}$. 
This is also a useful representation for mean-field limits. Indeed, when the number of agents tends to infinity, the concept of graph and neighbors is lost. Instead, the interaction potential, which can be based on the relative distance between agents, can be easily transported to the mean-field setting.

\subsection{Interaction Network in bounded-confidence Models } % Aylin, Sean
\input{Networks-BoundedConfidence}\label{Sec:Net_BoundConf}
\subsection{The Interaction Potential } \label{Sec:Net_Pot}

\input{Networks-InteractionPotential}

\section{Role of the State Space} 
In standard models, the agents evolve in Euclidean spaces $\R^{Nd}$ or $\R^{2Nd}$.
For modeling purposes, one might need to consider more complex state-spaces like compact manifolds, for instance $\Sp^1$ or $\T$. The dynamics then give rise to new kinds of equilibria that differ from the usual consensus or alignment. We will present such models for both first-order and second-order dynamics.

\subsection{First-order dynamics } % Aylin, Sean, Nastassia
\input{StateSpace-Opinion}

\subsection{Second-order dynamics } % Nastassia
\input{StateSpace-Animal}

\section{Control} % Marco
% Sparse control, stabilization, optimal control, emergence of consensus
\input{Control}

\section*{Acknowledgments}
The authors acknowledge the partial support of the NSF Project ``KI-Net", DMS Grant \# 1107444.

\bibliographystyle{abbrv}
%\nocite{*}
\bibliography{Bibliography.bib}

\end{document}

%% file: Intro.tex
% Introduction

A fascinating feature of large groups of autonomous agents is their ability to form organized global patterns even when individual agents interact only at a local scale. This is usually referred to as \textit{self-organization}. We use the term \textit{Social Dynamics} to indicate the study of such global behaviors, with an emphasis on  understanding the mechanisms leading from local rules to global phenomena, as well as identifying the resulting global pattern formation. 

Social dynamics models can be classified as first-order models and second-order models. In first-order models, we refer to the variables of interest as \textit{opinions}, even though such models can describe a wide range of attributes such as positions, market shares or wealth.   
The opinion of each agent is affected by neighboring agents' opinions in the state space. 
On the other hand, in second-order models, the variables of interest are the \textit{velocities}, obtained as the time derivatives of the positions. Each agent's velocity is affected by the velocities of agents whose positions are close in the state-space. 
%Second-order models are frequently used to describe animal group behavior such as flocking.

%Social dynamics encompasses many specific fields of application. In this review we will focus mostly on Opinion Dynamics, Animal Groups and Crowd Dynamics.   

First-order models (or opinion dynamics) 
can give rise to patterns such as
%consists of modeling opinion formation to understand trends such as 
\textit{consensus} (i.e. agreement of all states), \textit{polarization} (i.e. disagreement between two opposite parties) or \textit{clustering} (i.e. break-down of the opinions into several subsets).
%agreement (i.e. consensus), disagreement (polarization) or party formation (clustering). 
A first formulation of opinion dynamics can be traced back to French's research on social influence 
\cite{F56}, followed by works by 
Harary \cite{H59}, 
De Groot \cite{DG74} and 
Lehrer \cite{L75}, all focusing on linear models. 
More recently, nonlinear models were introduced and analyzed by Krause \cite{K97, K00}, Dittmer \cite{D00}, Hegselmann and Flache \cite{HF98}. 
 
Second-order models are commonly applied to animal groups to study coordinated collective behavior (as done by Couzin et al. \cite{C02}, Cristiani, Frasca and Piccoli \cite{CFP11}, Giardina \cite{G08}, Krause and Ruxton \cite{KR02}, Leonard \cite{L13} and Sumpter \cite{S06}) for example in fish (Huth and Vissel \cite{HW92}, Parrish, Viscido and Grunbaum \cite{P02}) or birds (Ballerini et al. and Cucker and Smale \cite{BCC08, CS07}). 
Some models have been designed to include simple interaction rules like attraction, short-distance repulsion and mimetic orientation or alignment. 
Agreement of all agents in the velocity variable is referred to as \textit{alignment} or \textit{flocking}. 

% \cite{B09, C02, G08, HW92, KR02, P02, S06}.

%(for instance representing animal groups) have also been extensively studied not only by mathematicians but also by biologists, ethologists, physicists and engineers. The main interest in this field is to understand the phenomenon of \textit{self-organization}, i.e. how simple interactions lead to complex animal behavior at the scale of a group, which led to the works of \cite{B09, C02, G08, HW92, KR02, P02, S06}. The basic interaction rules have been designed to comprise attraction terms, short-distance repulsion terms and mimetic orientation or alignment terms such as some second-order models focusing on the phenomenon of alignment \cite{CFPT13, CFPT15, CS07, FPR14, V95}.

%As a third emphasis, we will focus on the field of Crowd Dynamics, with some very concrete applications involving anticipation and control of crowd behavior \cite{HFV00, HM95, KM03}, and regulation of pedestrian flows in public environments \cite{ CPT14, CPT11}. Such work can have a direct impact on architecture, escape planning, panic control and event organization.  
%The local interaction rules of crowd dynamics are typically more complex than those of opinion dynamics or animal groups. For instance, Helbing and Molnar proposed a microscopic model based on ``social forces" that measure the internal motivations of the individuals \cite{HM95}. Self-organization also arises from such models, leading to interesting behaviors such as lane formation \cite{KM03}. 

The aim of this survey is to describe the role of two elements affecting the dynamics for both first-order and second-order models: the interaction network and the state space. We will also explore ways to control the dynamics to drive the system to a desired state.
% We will also pose and illustrate open problems related to these factors.  

The interaction network plays a critical role in the emergence of global patterns. Depending on the network, opinion formation models may lead to consensus among all the opinions or to the formation of separate clusters. 
The system's dynamics and the network's dynamics may be coupled.
For instance, bounded-confidence models allow agents to interact only if they are within a certain radius of each other in the state-space, as proposed by Hegselmann and Krause \cite{HK02}. On the other hand, it was shown that heterophilious dynamics enhances consensus by Motsch and Tadmor \cite{MT14}. Another distinction can be made between metric and topological interactions. A network based on metric interactions links agents based on their distance in the state space, whereas one based on topological interactions links an agent to another if it is among its $k$ closest neighbors, which can lead to asymmetric relations and interesting patterns. Furthermore, a network may be constant in time or time-dependent.
%In this case the network dynamics and systems ones are completely coupled.

The state-space is another factor that greatly influences the dynamics. Most studies have considered dynamics in Euclidean spaces (most often 1-dimensional for opinion models
and 2 or 3-dimensional for animal groups). One can also study the same dynamics on general Riemannian manifolds. For instance, a nonlinear model of opinion formation on the sphere was studied by Caponigro, Lai and Piccoli \cite{CLP15}, with a rich structure leading to unusual equilibria. These models are based on the projection of the linear dynamics in the ambient space onto the tangent space of the manifold. Consensus dynamics on special orthogonal groups were also studied, for example by Sarlette and Sepulchre \cite{sarlette2009, sepulchre2011}, motivated by applications to satellites or ground vehicles. 
%Furthermore, results depend on the choice of the manifold metric. 
%We will explore the influence of the metric, particularly the differences between the distance inherited from the metric of the ambient space and the natural metric of the manifold. 

%Lastly, external agents may act on the system to control it. 
A large number of applications involve control of robotic networks or autonomous vehicles, as done by Bullo, Cort\'es, and Mart\'inez \cite{B09}.
Control is used to impose consensus or alignment when it is not reached naturally (see Caponigro, Fornasier, Piccoli, and Tr{\'e}lat \cite{CFPT13, CFPT15}), or to guide the agents in a specific direction, as done by Leonard for the migration of animal groups \cite{L13}.
Ways of controlling the system include spreading leaders among the group or acting on the network. Due to the high dimensionality of Social Dynamics systems, control can be excessively demanding in computational resources. It is then convenient to consider the mean-field limit of the system. Numerous theories have been developed to control the resulting kinetic equation. Some approaches require taking the limit (in some sense) of the finite-dimensional controlled system. For example, Fornasier and Solombrino have introduced a concept of $\Gamma$-limit for optimal control problems \cite{FS14}, and Fornasier, Piccoli and Rossi have extended the idea of control by leaders \cite{FPR14}. One can also control the PDE directly, as done by Piccoli, Rossi and Tr\'elat \cite{PRT15}. Other approaches involve controling the interaction kernel (see Albi, Herty and Pareschi \cite{AHP15}), or using mean-field games, a theory developed by Lasry and Lion \cite{LL07} and Caines \cite{C13}.

%% file: Overview.tex
In this section we give general definitions of the concepts that we will use. We also provide some examples of common first-order and second-order systems and distinguish between finite-dimensional and infinite-dimensional models.

\subsection{General Notations and Definitions}\label{Sec:Overview-gen}

In the following chapter we will differentiate between two branches of models: 
\begin{itemize}
\item First-order models (also referred to as opinion dynamics) that can lead to \textit{consensus}
\item Second-order models (mostly related to animal group models) that can lead to \textit{flocking} or \textit{alignment}
\end{itemize}

We shall write first-order dynamics as follows: 
\begin{equation}\label{eq:syst-gen}
\dot{x}_i = \frac1N \sum_{j\in\NN_i} a_{ij}\; (x_j-x_i), 
\quad i\in\{1,...,N\},
\end{equation}
and second-order systems as follows:
\begin{equation}\label{eq:syst2-gen}
\begin{cases}
\dot{x}_i = v_i \\
\dot{v}_i = \frac1N \sum\limits_{j\in\NN_i} a_{ij}\;(v_j-v_i)
\end{cases}
\quad i\in\{1,...,N\},
\end{equation}
where $N$ is the number of agents, $x_i\in\R^d$ is the position of agent $i$ in the state space, $v_i\in\R^d$ is its velocity,
$\NN_i$ is the set of agents interacting with agent $i$ and 
$a_{ij}$ are interaction coefficients for each pair of agents $(i,j)$. They form the interaction matrix $A = (a_{ij})_{i,j\in\{1,...,N\}}$.
%:=a(\|x_i-x_j\|)$ is a function of the distance between the agents $i$ and $j$ in the state space (for metric interaction networks). 
Unless otherwise specified, we consider that first-order systems evolve in $\R^{Nd}$ (where $d$ is the dimension of the state-space) and second-order systems are in $\R^{2Nd}$.
%A further distinction must be made in the type of interactions. 
\begin{remark}
The dynamics \eqref{eq:syst-gen} and \eqref{eq:syst2-gen} can be written in a more general form: $\dot{x}_i = \frac{1}{\deg_i} \sum_{j\in\NN_i} a_{ij}\; (x_j-x_i) $  or 
$\dot{x}_i = v_i ; \quad
\dot{v}_i = \frac{1}{\deg_i} \sum_{j\in\NN_i} a_{ij}\;(v_j-v_i)$
 where $\deg_i$ is a scaling factor. Typical choices for scaling factors are: $\deg_i=N$,  $\deg_i=\mathbf{card}(\NN_i)$ 
 or  $\deg_i=\sum_{j\in\NN_i} a_{ij}$, where $\mathbf{card}(\cdot )$ denotes the cardinality of a set.
\end{remark}

First-order models are also referred to as \textit{consensus} models. Consensus is an equilibrium state in which all agents have the same opinion: $x_i=x_j$ for all $i, j\in\{1,...,N\}$.
Second-order models are also called \textit{alignment} (or \textit{flocking}) models. Alignment or flocking is the equilibrium set in which all agents have the same velocity: $v_i=v_j$ for all $i, j\in\{1,...,N\}$. For this reason, the velocity $v$ is also referred to as the \textit{consensus variable}, to distinguish from the position $x$.

%Two approaches can be distinguished. If the sets of interacting agents $\NN_i$ are strict subsets of the total set of agents $\{1,..., N\}$, then 
The system can be viewed as a network represented by a (possibly time-varying) directed weighted graph $\GG = (\VV,\EE)$. We define the set of vertices $\VV=(\nu_i)_{i\in\{1,...,N\}}$ corresponding to the set of agents, and the set of edges $\EE \subseteq \VV\times\VV$, so that an edge exists between two vertices $i$ and $j$ if and only if $a_{ij}\neq 0$. The edges are weighted by the interaction coefficients $a_{ij}$. 

Most often, the interaction coefficients are defined by an interaction potential $a(\cdot)$ such that $a_{ij}:=a(\|x_i-x_j\|)$. When modeled as such, the strength of interaction is a function of the distance between agents in the position space. This generates a fundamental difference between first-order and second-order models. In first-order models, the variable of interest is the position and its tendency to agree with other agents' positions depends on the distance between the agents. In second-order models, the velocity's tendency to align with other agent's velocities depends on the difference in their positions. 

If the interactions between agents are only local, that is if an agent interacts exclusively with close neighbors, we refer to \textit{bounded confidence} models, a term first introduced by Hegselmann and Krause \cite{HK02}.
Then $\NN_i$ denotes the set of closest neighbors of the $i$-th agent. Bounded confidence models will be examined in Section \ref{Sec:Net_BoundConf}.
We look in particular at two ways to define proximity of agents. 
In the case of bounded confidence with metric interaction, given a radius $r>0$,
\begin{equation}\label{eq:Nimetric}
\NN^r_i(x) = \{j\in\{1,...,N\},  \|x_i-x_j\|\leq r\}
\end{equation}
In the case of bounded confidence with topological interaction, we define the relative separation between two agents as
$\alpha_{ij}=\mathbf{card} \{k:\|x_i-x_k\|\leq\|x_i-x_j\| \}$.
 Then the set of neighbors of agent $i$ is defined as the set of its $k$ closest neighbors, i.e. 
 \begin{equation}\label{eq:Nitopo}
 \NN^k_i(x) = \{j\in\{1,...,N\}, \alpha_{ij}\leq k\}, 
 \end{equation}
 for a given $k\in \N$. 

When all agents interact with all others,
$\NN_i=\{1,...,N\}$ for all $i\in\{1,...,N\}$. Then the network is fully connected but its edges may have varying weights. The behavior of the system will depend on the interaction potential $a(\cdot)$, as seen in Section \ref{Sec:Net_Pot}.

% systems \eqref{eq:syst-gen} and \eqref{eq:syst2-gen} can be rewritten as:
%\begin{equation}\label{eq:syst}
%\dot{x}_i = \frac1N \sum_{j=1}^{N} a_{ij}\;(x_j-x_i), 
%\quad i\in\{1,...,N\},
%\end{equation}
%and  
%\begin{equation}\label{eq:syst2}
%\begin{cases}
%\dot{x}_i = v_i \\
%\dot{v}_i = \frac1N \sum\limits_{j=1}^{N} a_{ij}\;(v_j-v_i)
%\end{cases}
%\quad i\in\{1,...,N\}.
%\end{equation}

%In the case of non-Euclidean state-spaces, the $d$-dimensional sphere will be denoted $\Sp^d$, the torus $\T^d$ and a general manifold $\MM$. 

%Notations for control of the system will be taken from \cite{CFPT13}. For instance, an additive control in \eqref{eq:syst2} will be written as: 
%\begin{equation}\label{eq:syst2-cont}
%\begin{cases}
%\dot{x}_i = v_i \\
%\dot{v}_i = \frac1N \sum\limits_{j=1}^{N} a_{ij}\;(v_j-v_i) + u_i
%\end{cases}
%\quad i\in\{1,...,N\},
%\end{equation}
%with possible constraint on the control: $\sum_{i=1}^N \|u_i\| \leq M$ for some constant $M\in\R$. 

\subsection{Examples of first-order consensus models}\label{Sec:Overview-order1}

We start by giving two examples of common first-order consensus models. The \textit{Voter} model is a discrete-time model, whereas the Hegselmann-Krause model is a system of ODEs.

\subsubsection{The Voter model}
One system used to explore the dynamics of cellular automata is the Sznajd Model (SM) \cite{behera}. The SM is an example of discrete-time and discrete-state model. 
The alignment variable of each agent can take one of two values, referred to as \textit{spin up} or \textit{spin down}.
%, and some agents remain ``on the fence'' as they can have a neutral opinion. 
The dynamics of this particular model operate on a one or two-dimensional lattice.  In this system, the agents change their opinion (spin up or spin down) based on specific interaction rules: the \textit{ferromagnetic} interaction (that is, if $x_{i} = x_{i+1}$ then at the next step adjacent agents will satisfy with a given probability $x_{i-1} = x_{i} = x_{i+1}= x_{i+2}$) and the \textit{antiferromagnetic} interaction 
 (if $x_{i} = - x_{i+1}$ then an antisymmetric pattern forms: $- x_{i-1} = x_{i} = - x_{i+1}= x_{i+2}$). The model has been extended to higher dimensional opinion and complex network topologies.  The motivation for this model comes form the postulate that ``agreement  generates agreement'', that is, if two agents reach a consensus then all agents directly connected to them are induced to agree. In other words, in the Sznajd model, the opinion flows out from a group of agreeing agents, a concept known as \textit{social validation}. 

%   \begin{enumerate}
%   \item  A $2\times2$ panel of four agents, if not all four spins are parallel, will leave its six neighbors unchanged.
%   \item A neighboring pair will persuade the pair's six neighbors to adopt the same spin.
% \end{enumerate}
%  The system will always reach a consensus with these interaction rules.  The simple dynamics are of interest because they capture some aspect of conformity in social networks.  Neighbors will adopt behavior, but only if the behavior is sufficiently and locally prevalent.
  In \cite{behera}, the authors show that the SM is a special case of a linear Voter Model.  The Voter Model (VM) is one of the simplest mathematical models of cooperative behavior, and its dynamics are well understood.  Here, each node of a graph begins as either one of two states: spin up or spin down.  
 The system then follows an algorithm:
 \begin{enumerate}[nolistsep]
  \item  pick a random voter 
   \item  the selected voter adopts the state of a randomly chosen neighbor
   \item repeat steps 1 and 2 until consensus
 \end{enumerate}
 Once the system reaches consensus, all nodes are spin up or spin down. In this system, the interactions between a randomly selected voter $x_i$ and a randomly chosen neighbor $x_j$ is $x_i = x_j$.  In other words, the interaction is described by complete agreement with one of the neighboring agents.

\subsubsection{The Hegselmann-Krause model}

The Hegselmann-Krause model (HK) is a classical example of a first-order nonlinear opinion formation model \cite{HK02}. 
It was designed in the context of opinion dynamics, 
and captures well-known phenomena
such as formation of consensus and emergence of clustering. 
Agents modify their own opinion to average neighboring opinions as follows:

\begin{equation}\label{eq:HK}
\dot{x}_i = \frac{1}{\mathbf{card}(\NN_i)} \sum\limits_{j\in\NN_i}(x_j-x_i) \quad \text{ for all } i\in\{1,...,N\}, \quad x_i\in\R^d,
\end{equation}
where 
$\NN_i = \{j : \|x_i - x_j\|\leq r\}, \; r > 0$, is the set of agents interacting with agent i. The radius $r$ can be interpreted as the level of confidence. This model captures the fact that an individual tends to trust only opinions that do not differ from its own by more than $r$.
Since the interaction region is bounded, the HK
model is also called \textit{bounded confidence} model. Depending on the size of the interaction
regions and the density of agents in the domain, different phenomena are observed. If the
interaction is strong enough (i.e. $r$ is big enough), the agents can be brought to consensus,
i.e. convergence to a single opinion. If the interaction regions are too restricted, one observes
clustering around different opinions. 
A wide variety of models have been developed by varying the confidence region $\NN_i$. Hegselmann and Krause have for instance looked at (one-dimensional) asymmetric confidence: $\NN_i = \{j : -r_l\leq x_i - x_j \leq r_r \}, \; r_l > 0, r_r>0$ \cite{HK02}.
Recently, Motsch and Tadmor have analyzed models with interaction strength increasing with the distance between agents, showing that this so-called heterophilious dynamics enhances consensus \cite{MT14}.

Jabin and Motsch have studied a slightly different model for opinion formation, that can be written as: 
\begin{equation}
\dot{x}_i = \frac{ \sum_{j} \phi_{ij}\; (x_j-x_i)}{\sum_{j} \phi_{ij}}, 
\quad i\in\{1,...,N\}, \quad x_i\in\R^d,
\end{equation}
where $\phi$ is the influence function and we define $\phi_{ij}:=\phi(\|x_i-x_j\|^2)$ \cite{JM14}. One can prove that under appropriate conditions on the influence function, the system leads to clustering. For instance, if 
\begin{itemize}
\item $\phi\in L^\infty(\R^d)$ with compact support in $[0,1]$ 
\item for any $\epsilon>0$, $\phi\in W^{1,\infty}([0,1-\epsilon])$ and $\phi$ is strictly positive on $[0,1-\epsilon]$
\item $|\phi'(r)|^2\leq C\phi(r)$ for all $r\in [0,1]$
\end{itemize}
then there exists a set of cluster centers $\{\bar{x}_i\}$ such that for all $i$, $x_i(t)\rightarrow_{t\rightarrow\infty}\bar{x}_i$, and for any $i,j$, either $\bar{x}_i=\bar{x}_j$ or $|\bar{x}_i-\bar{x}_j |\geq 1$ \cite{JM14}.

In one dimension, we can even characterize the rate of convergence to the clusters in the following way. Assume that $d=1$ and $\phi\in W^{1,\infty}([0,1))$ with $\inf_{[0,1)}\phi>0$. Then for each agent $i$ there exists $\bar{x}_i$ depending on the initial positions of all the agents such that 
$|x_i(t)-\bar{x}_i|\leq C e^{-\lambda(t-t_0)}$ for all $t \geq t_0$, 
where the constants $C$ and $\lambda$ are determined a priori by the total number of agents $N$ and by the influence function $\phi$, and the time $t_0$ depends on $N$, $\phi$ and the diameter of the initial support \cite{JM14}.

\subsection{Examples of finite-dimensional and infinite-dimensional second-order alignment models}\label{Sec:Overview-order2}

There exists a wide variety of second-order models, that have been developed mainly to describe the behaviors of animal groups or robotic networks. Some early models like the Vicsek model \cite{VCABC95} are defined in discrete time, and require to update each agent's state at successive time intervals. Other models like the Cucker-Smale one \cite{CS07} are continuous in time and require the use of ODE's. 
We also look at the limit of such models when the number of agents tends to infinity, which is referred to as the \textit{mean-field limit}. 

\subsubsection{The Vicsek model}
A classic example of discrete-time model is the Vicsek model \cite{VCABC95}, proposed to describe interactions within animal groups such as a school of fish. It represents each agent (or fish) by its position $x_k$ and its velocity angle $\theta_k$, all velocities having constant norm $v$. The positions and angles are updated in the following way: 
\begin{equation}
\begin{cases}
x_k(t+\Delta t) = x_k(t) + v_k(t)\Delta t \\
\theta_k(t+\Delta t) = \langle\theta(t)\rangle_r + \Delta \theta_k
\end{cases}
\end{equation}
where $\Delta \theta_k$ is a noise term, and $\langle\theta(t)\rangle_r$ represents the average direction of the velocities of particles being within a circle of radius $r$ of particle $k$.\\

\subsubsection{The Finite-dimensional Cucker-Smale model} 
The prototypical second-order model for the interaction of $N$ agents is the \textit{Cucker-Smale model} (CS) \cite{CS07}:
\begin{equation}\label{eq:CS}
\begin{cases}
\dot{x}_i(t)&=v_i(t)\\
\dot{v}_i(t)&=\frac{1} {N}\sum\limits_{j=1}^N a(\Vert x_{j}(t) - x_{i}(t)\Vert )(v_j(t)-v_i(t)),
\end{cases}
\qquad\qquad i=1,\ldots,N
\end{equation}
where  $x_i\in\mathbb{R}^d$, $v_i\in\mathbb{R}^d$, and $a \in C^{1}([0,+\infty))$ is a nonincreasing positive function called interaction potential or rate of communication. 
 In the classical CS model, we have $a(s)=\frac{1}{(1+s^2)^\beta}$, with $\beta>0$.
Here, $x_i$ is the main state of the agent $i$, and $v_i$ is its consensus parameter.
This model was initially introduced to describe the formation and evolution of languages, and was then also used for describing the flocking of a swarm of birds \cite{CS07} or spacecraft formation \cite{PEG09}.

We consider the finite-dimensional CS model \eqref{eq:CS} and provide some results published by Caponigro, Fornasier, Piccoli and Tr\'elat \cite{CFPT13}.
We define the space barycenter $\bar x(t)$ and the mean velocity $\bar v(t)$ by
$$
\bar x(t) = \frac{1}{N} \sum_{i=1}^{N} x_{i}(t),\qquad
\bar v = \frac{1}{N} \sum_{i=1}^{N} v_{i}(t).
$$
Then $\dot {\bar x}(t)={\bar v}$ and the mean velocity is constant: ${\bar v}(t) = {\bar v}(0)$ for all $t$.
Define the spatial variance by
$$
X(t) =  \frac{1}{2N^{2}} \sum_{i,j=1}^N \Vert x_{i}(t) - x_{j}(t)\Vert ^{2} .
$$
The velocity variance is
\begin{equation}\label{eq:CSvelvariance}
V(t) = \frac{1}{2N^{2}} \sum_{i,j=1}^N \Vert v_{i}(t) - v_{j}(t)\Vert ^{2} 
= \frac{1}{N} \sum_{i=1}^{N}\Vert v_i(t)-\bar v\Vert ^{2} ,
\end{equation}
and we have:
$$
\dot{V}(t) = - \frac{1}{N} \sum_{i,j=1}^N  a(\Vert x_{j}(t) - x_{i}(t)\Vert ) \Vert v_i(t)-v_j(t)\Vert^2   \leq 0.
$$

\begin{definition}
A solution $(x(t), v(t))$ converges to \emph{alignment} (or \textit{flocking}) if
\begin{itemize}
\item[(i)] there exists $X_M>0$ such that $X(t)\leq X_M$ for every $t>0$,
\item[(ii)] $v_{i}(t) \underset{t \to +\infty}{\longrightarrow} \bar v$, for every $i=1,\ldots,N$, or equivalently,  
$V(t) \underset{t \to +\infty}{\longrightarrow}  0$.
\end{itemize}
\end{definition}

Note that, since $a$ is nonincreasing, we have $\dot{V}(t) \leq -2a\left(\sqrt{2N X(t)}\right) V(t)$, and hence \textit{if $X(t)$ remains bounded} then $\dot V\leq -\alpha V$, which implies flocking. But the difficulty is that the group of agents does not necessarily remain confined, and for this reason convergence to alignement is not guaranteed. More precisely, Ha, Ha and Kim provided the following result \cite{HHK10}:

\begin{prop} \cite{HHK10}  \label{propHaHaKim}
Let $(x_{0},v_{0}) \in (\mathbb{R}^d)^N \times (\mathbb{R}^d)^N$ be such that
$$\sqrt{V(0)} \leq \int_{\sqrt{X(0)}}^{+\infty} a (\sqrt{2N}r) \, dr .$$
Then the solution with initial data $(x_{0},v_{0})$ tends to alignment.
\end{prop}

\begin{figure}[h]
\begin{center}
\begin{tikzpicture}
\draw[scale = 0.5, thick,->] (0,0) -- (11,0) node[right] {$X_0$};
\draw[scale = 0.5, thick,->] (0,-0.1) -- (0,11) node[left] {$V_0$};
%\draw[scale = 1,domain=0:1,smooth,variable=\x,blue] plot ({\x},{ 1/(200)*
%((3.14)/2- atan(sqrt(\x)))^2 });
\draw[scale = 0.5,domain=0.01:10,smooth,variable=\x,black] plot[id=x,samples=200] ({\x},{4*((3.14)/2- rad(atan(sqrt(\x))))^2});
\fill [scale = 0.5, gray, domain=0:10, variable=\x,samples=200]
      (0, 0)
      -- plot ({\x}, {4*((3.14)/2- rad(atan(sqrt(\x))))^2})
      -- (10, 0)
      -- cycle;
\draw	(0,-0.1) node[above right] {\color{white} \emph{Flocking region}};
\draw	(0,0) node[below left] {0};
\end{tikzpicture}
\caption{\small Flocking region $\sqrt{V_0}\leq\frac{1}{\sqrt{2N}}(\frac{\pi}{2}-\arctan(\sqrt{X_0}))$ corresponding to the CS model \eqref{eq:CS} with parameter $\beta=1$, see Proposition \ref{propHaHaKim}.}\label{fig_consensusregion}
\end{center}
\end{figure}
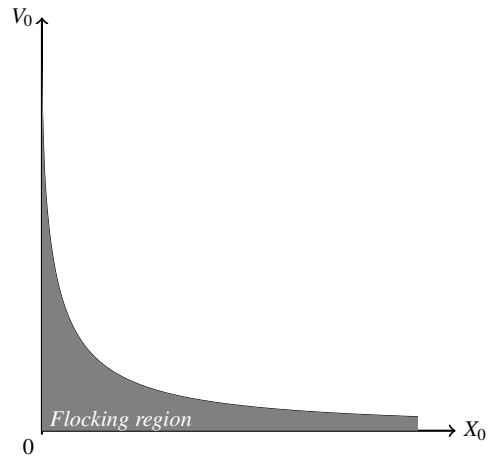

Figure \ref{fig_consensusregion} illustrates the \textit{self-organization of the group} in what we can call the \textit{flocking region} (region of natural asymptotic stability to flocking).

\subsubsection{The Infinite-dimensional kinetic Cucker-Smale model}\label{sec_kineticCS}
In numerous applications such as risk-taking in economics, pricing models and opinion formation, the system is made of a very large number of agents.
Studying and simulating social dynamics systems becomes a particularly challenging problem when the dimension of the system increases. This is referred to as \textit{the curse of dimensionality}, a term coined by Bellman in the context of dynamic optimization of high-dimensional systems.
One way around this problem is to move away from the microscopic viewpoint where each agent is considered individually, and consider instead the mean-field limit, which provides a kinetic description of the system. This approach consists of 
approximating the influence of all agents on any given individual by one averaged effect.
% When the number of agents is large, it is convenient to study the distribution of the particles positioned at $(x,t)$ with velocity $v$, instead of the motion of each individual particle.
  Derivation of kinetic models have been intensively studied, for example by 
  Ca\~{n}izo, Carillo and Rosado \cite{CCR11}, by Ha and Tadmor for the CS model \cite{HT08} or by Degond and Motsch for the Viscek model \cite{DM08_2, DM08}.
  
 When the number of agents is large, one often refers to the agents as \textit{particles}.
 Let $\mu(t,x,v)$ denote the distribution function of particles positioned at $x\in\R^d$ at time $t>0$ with velocity $v\in\R^d$.
By taking the mean-field limit in system \eqref{eq:CS}, we obtain the \textit{kinetic Cucker-Smale} model:
\begin{equation}\label{kineticCS}
\partial_t \mu +\langle v, \nabla_x \mu\rangle+\mathrm{div}_v \left( \xi[\mu]\, \mu \right)=0
\end{equation}
where $\mu(t)$ is a probability measure on $\mathbb{R}^d\times\mathbb{R}^d$ 
 (if $\mu(t,x,v)=f(t,x,v)\, dx\, dv$, then $f$ is the density of the particles), 
and $\xi[\mu]$ is the \textit{interaction kernel}, given by
$$
\xi[\mu](x,v)=\int_{\mathbb{R}^d\times\mathbb{R}^d} a(\Vert x-y\Vert) (w-v)\,d\mu(y,w).
$$
The link with the finite dimensional system is given by the empirical measure
$$
\mu(t)=\frac{1}{N}\sum_{i=1}^N\delta_{(x_i(t),v_i(t))}.
$$
Indeed, plugging this measure in \eqref{kineticCS}, we find that $(x_i(t),v_i(t))$  satisfy exactly \eqref{eq:CS}.
The kinetic equation \eqref{kineticCS} can be written as
$$
\partial_t \mu + \mathrm{div}_{(x,v)} \left( V[\mu]  \mu \right) = 0 ,
$$
with the velocity field
$$
V[\mu] = \begin{pmatrix} v \\ \xi[\mu] \end{pmatrix} .
$$
The so-called \textit{particle flow} $\Phi(t)$ generated by $V[\mu(t)]$ yields the \textit{characteristics}
$$
\dot x(t) = v(t),\qquad \dot v(t) = \xi[\mu(t)](x(t),v(t)).
$$
The motion of any such particle follows exactly the finite-dimensional CS system. This justifies the wording \textit{particle}. Moreover, the solution to the kinetic equation \eqref{kineticCS} is formally:
$$
\mu(t)=\Phi(t)\#\mu^0 ,
$$
that is, the pushforward under the flow $\Phi(t)$ of the initial measure.

Similarly to the finite-dimensional case, we present some of the properties of the infinite-dimensional model. 
In the infinite-dimensional setting, we define the space barycenter and mean velocity by
$$
{\bar x}(t)=\int_{\R^d\times\R^d} x \, d\mu(t)(x,v),\qquad
{\bar v}=\int_{\R^d\times\R^d} v \, d\mu(t)(x,v) .
$$
Then $\dot {\bar x}(t)={\bar v}$ and $ \dot{\bar v}= 0$, as in finite dimension.
Defining (as before) the spatial and velocity variances by
$$
X(t) = \int_{\R^d\times\R^d} \Vert x-{\bar x}(t)\Vert^2\,d\mu(t)(x,v),\qquad
V(t)=\int_{\R^d\times\R^d} \Vert v-{\bar v}\Vert^2\,d\mu(t)(x,v)  ,
$$
we have
$$
\dot{V}(t) = - \iint  a(\Vert x - y\Vert ) \Vert v-w\Vert^2\, d\mu(t)(x,v)\, d\mu(t)(y,w) \leq 0.
$$
We expect that $V(t) \underset{t\rightarrow+\infty}{\longrightarrow} 0$, but as in finite dimension, this is not guaranteed, unless the population of agents remains confined. This justifies the following definition. The notation $\mathrm{supp}$ stands for the support of a measure.

\begin{definition}
A solution $\mu\in C^0(\R,\mathcal{P}_c(\R^d\times\R^d))$ converges to \textit{alignment} (or \textit{flocking}) if:
\begin{itemize}
\item[(i)] there exists $X_M>0$ such that $\mathrm{supp}(\mu(t))\subseteq B({\bar x}(t),X_M)\times \R^d $ for every $t>0$,
\item[(ii)] $V(t) \underset{t\rightarrow+\infty}{\longrightarrow} 0$.
\end{itemize}
\end{definition}

Piccoli, Rossi and Tr\'elat \cite{PRT15} provided the infinite-dimensional counterpart of Proposition \ref{propHaHaKim}. As in finite dimension, it defines a consensus region, that is, a set of initial conditions for which the group of agents will naturally converge to alignment:

\begin{prop} \cite{PRT15}
Let $\mu^0\in \mathcal{P}_c(\R^d\times\R^d)$. Define the space and velocity barycenters  ${\bar x}^0=\int x\,d\mu^0,\ \ {\bar v}=\int v\,d\mu^0$ and the space and velocity support radii:
\begin{equation*}
\begin{split}
X^0=&\inf\left\{X\geq 0\mid \mathrm{supp}(\mu^0)\subset  B({\bar x}^0,X)\times \R^d\right\},\quad \\
V^0=&\inf\left\{V\geq 0\mid \mathrm{supp}(\mu^0)\subset  \R^d\times B({\bar v},V)\right\}.
\end{split}
\end{equation*}
If
$$
V^0< \int_{X^0}^{+\infty} a(2 x)\,dx,
$$
then the solution $\mu(t)$ with initial datum $\mu(0)=\mu^0$ converges to consensus.
\end{prop}

%% file: Networks-BoundedConfidence.tex
In this section we study the influence of the interaction network in bounded-confidence models. We review known properties of such models, propose open problems concerning the equilibrium sets, and provide numerical simulations illustrating the known and conjectured properties.

\subsubsection{Properties of bounded-confidence models}

The rationale for bounded confidence models is that it is unlikely for one agent to be influenced by another one whose opinion is too far from its own.  This kind of interaction gives rise to clusters of opinions (see for instance~\cite{BHT}). We also mention the bounded confidence model by Deffuant, see~\cite{deffuant} in which the opinions belong to real intervals too but the pairs of interacting agents are chosen randomly.

\begin{figure}[h!]
\begin{center}
           \subfloat[Metric ]{\includegraphics[width = 0.24 \textwidth, height=0.20 \textwidth]{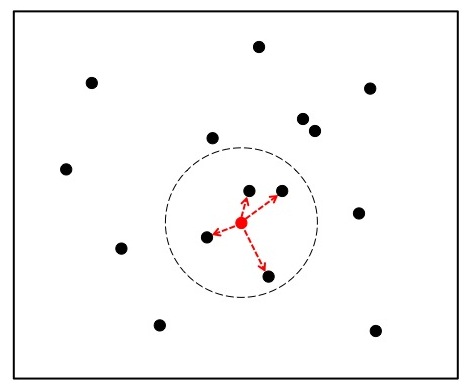}\label{fig:Nimetric}}\hskip 0.3in 
           \subfloat[Topological]{\includegraphics[width = 0.24 \textwidth, height=0.20 \textwidth]{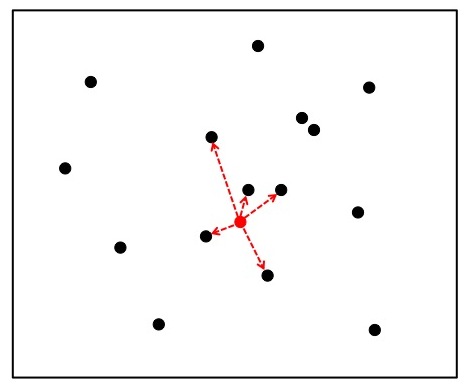}\label{fig:Nitopo}}\hskip 0.3in 
           \subfloat[Long-distance]{\includegraphics[width = 0.24 \textwidth, height=0.20 \textwidth]{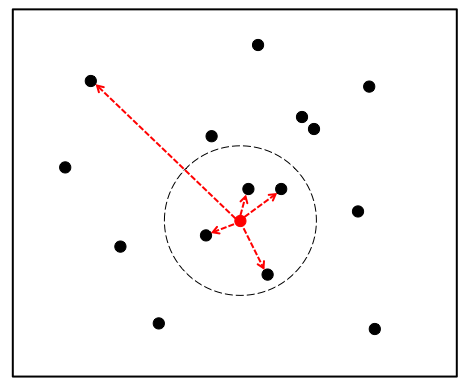}\label{fig:Nilong}}\hskip 0.3in 
        \caption{\small Representation of interacting neighbors for one agent according to the different interaction networks.  In (c), the long-distance connection is added to metric local interactions.}
\label{fig:Ni}
\end{center}
\end{figure}

As mentioned in Section \ref{Sec:Overview-gen}, two main types of interaction networks have been proposed in the literature. 
 In metric interaction networks, agents interact depending on their distance in the state space \cite{H13}: given a confidence radius $r>0$, we can define the interaction neighborhood $\NN_i^r$ \eqref{eq:Nimetric} , see Fig. \ref{fig:Nimetric}.
 In topological interaction networks, agents interact depending on their relative separation.
 %, i.e. on the number of agents separating them.
 Given $k\in\N$, we can define the interaction neighborhood $\NN_i^k$ \eqref{eq:Nitopo}, see Fig \ref{fig:Nitopo} .

%Both types of networks are time-varying, which highlights the concept of dynamic graphs. 
%In particular, for the Hegselman-Krause dynamics \eqref{eq:syst-gen} and for the Cucker-Smale second-order dynamics \eqref{eq:syst2-gen}, one can define the set of the $i$-th agent's neighbors as $\NN^r_i$ as follows: 
%\begin{itemize}
%\item For metric interaction: $\NN^r_i=\{ j\in\{1,...,N\},\; \|x_i-x_j\|\leq r\}$. 
%\item For topological interaction: $\NN^k_i=\{ j\in\{1,...,N\},\; \#\{ m\in\{1,...,N\} \; | \; \|x_i-x_m\|\leq \|x_i-x_j\| \}\leq k\}$.
%\end{itemize}

Both topological and metric interactions are local interactions. %Kleinberg worked on the \textit{small world} problem, showing that 
Adding long-distance connections to local ones greatly reduces the network's diameter and facilitates the spread of information \cite{K00}.
This is justified by the ubiquitous idea that
social networks are of small diameter, a property also known
as the six degrees of separation or \textit{small-world effect} \cite{WS98}.
In particular, Kleinberg (see \cite{KJ00, K00}) showed
that single long-distance random connections in locally organized networks lead to efficient routing procedures for spreading information. 
%  An effort to model the so called ``Small World Phenomena'' \cite{KJ00}, implements short range interactions, as in a bounded confidence model, but also a long range interaction decided at random. 
The small world phenomenon is characterized by short paths (relative to the size of the network) connecting any two nodes in the network, as illustrated in Fig. \ref{fig:Nilong}. The model as presented in \cite{KJ00} describes nodes on a square lattice which interact with the four adjacent nodes in the lattice, as well as one long range interaction that randomly forms an edge between a node and another non-neighboring node with a probability proportional to $\rho^{-a}$, where $\rho$ is the Manhattan distance between the two nodes.  
%Social networks often display characteristics of the Small World Phenomena.  One way to incorporate this behavior in a model is to combine many short range interactions with few long range ones. 

\begin{figure}[h!]
\begin{center}
           \subfloat[Metric, N=2]{\includegraphics[width = 0.20 \textwidth, height=0.20 \textwidth]{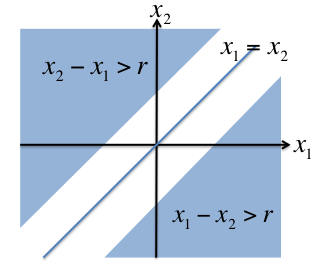}\label{fig:metricN2}}\hskip 0.3in 
           \subfloat[Top., k=2, N=5]{\includegraphics[width = 0.20 \textwidth, height=0.20 \textwidth]{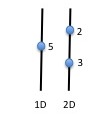}\label{fig:topoN5}}\hskip 0.3in 
           \subfloat[Top., k=2, N=5]{\includegraphics[width = 0.30 \textwidth, height=0.20 \textwidth]{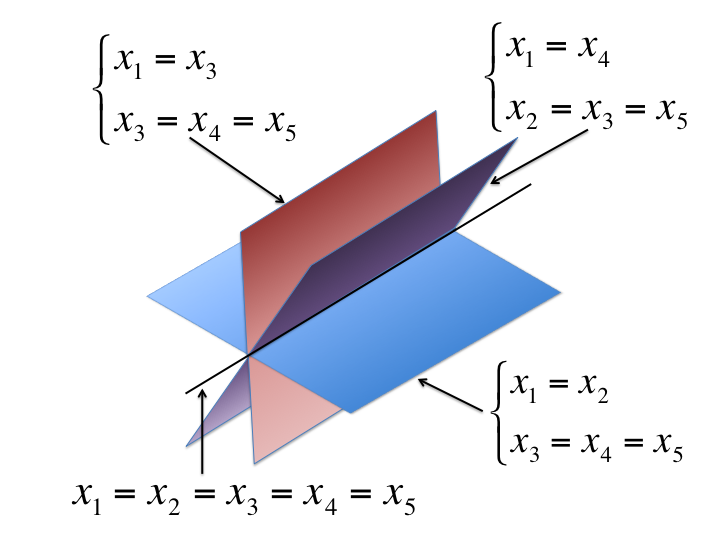}\label{fig:topoN5strata}}\hskip 0.3in 
          % \subfloat[Top., k=3, N=7]{\includegraphics[width = 0.20 \textwidth, height=0.20 \textwidth]{fig5.jpg}} \label{fig:topoN7}
        \caption{\small Equilibria for the HK system with metric and topological interactions for $d=1$. Figure (a) shows the equilibrium set in the metric case ($N=2$), with separate strata. Figure (b) shows the possible configurations for the agents' positions at equilibrium for the topological interaction ($k=2$, $N=5$), indicating the number of agents in each cluster and the dimension of the manifold. Figure (c) shows some of the non-separate strata of this equilibrium set.}
\label{fig:equilibria}
\end{center}
\end{figure}

\noindent\textbf{Equilibrium sets.}
To understand the mechanisms behind cluster formation, we studied equilibria for the HK dynamics, both with
metric and topological interactions.\\
\noindent
Let us start with metric interaction, with 2 or 3 agents in $\R$:
\begin{itemize}[nolistsep]
\item For $N=2$, the equilibrium set $E$ consists of 3 subsets: The line $x_1=x_2$; the half-plane $x_1-x_2>r$; the half-plane $x_2-x_1>r$ (see Figure \ref{fig:metricN2}).
\item In the case $N=3$, 13 equilibrium subsets can be enumerated: the line $x_1=x_2=x_3$; the 3 half-planes $\{x_i=x_j, \; x_k>x_i+r\}$; the 3 half-planes $\{x_i=x_j, \; x_k<x_i-r\}$ ; the 6 3D manifolds $\{x_i+r<x_j<x_k-r\}$
(with $i, j, k$ pairwise distinct in $\{1, 2, 3\}$).
\end{itemize}
Notice that in both cases, the equilibrium set is composed of pairwise disjoint manifolds with no common boundaries. We recall the following:
\begin{definition}
A set $E\subset\R^n$ is called \emph{stratified} in the sense of Whitney
\cite{W55}
if there exists a countable (locally finite) collection
of pairwise disjoint manifolds $M_i$, $i\in\N$, such that the following
holds:
\begin{enumerate}[nolistsep]
\item $M_i$ is an embedded manifold of dimension $d_i$.
\item If $M_i\cap\partial M_j\not=\emptyset$ then
$M_i\subset\partial M_j$ and $d_i<d_j$.
\end{enumerate}
Moreover we say that $E$ has separate strata
if for every $i\not=j$ we have $M_i\cap\partial M_j=\emptyset$.
\end{definition} 
We propose a general property for the equilibrium set:
\begin{conj}
For the HK dynamics \eqref{eq:HK}
with metric interaction, for all $d\in\N$ and $N\in\N$, the set of equilibria is a stratified manifold with separate strata.
\end{conj}

In the topological case, the number and nature of equilibrium sets depend on $k$. \\
If $k=1$ (i.e. there is no interaction between agents), the equilibrium set is $\R^N$ itself. 
If $k=2$ (each agent interacts with one other), we have to distinguish cases:
\begin{itemize}[nolistsep]
\item for $N=2$ or $N=3$, the equilibrium sets are respectively the lines $x_1=x_2$ and $x_1=x_2=x_3$. 
\item for $N\geq 4$, the equilibrium sets are more complex as they are composed of several manifolds. For instance, in the case $N=5$, the equilibrium set consists of the line $x_1=x_2=x_3=x_4=x_5$ and the ${5 \choose 2}=10$ half planes $\{x_i=x_j; \; x_k=x_l=x_m\}$ with $i,j,k,l,m$ pairwise distinct in $\{1,...,5\}$ (Fig. \ref{fig:topoN5}, \ref{fig:topoN5strata}). Notice that the line is in the boundary of all half planes.
\end{itemize}
Hence we propose the following:
\begin{conj}
For the HK dynamics \eqref{eq:HK} with topological interaction, for any $d\geq 2$ and $N\geq 4$, the set of equilibria is a stratified manifold 
with non-separate strata.
\end{conj}

\noindent

\subsubsection{Numerical results}

\begin{figure}[h!]
\begin{center}
           \subfloat[Metric ]{\includegraphics[width = 0.26 \textwidth, height=0.24 \textwidth]{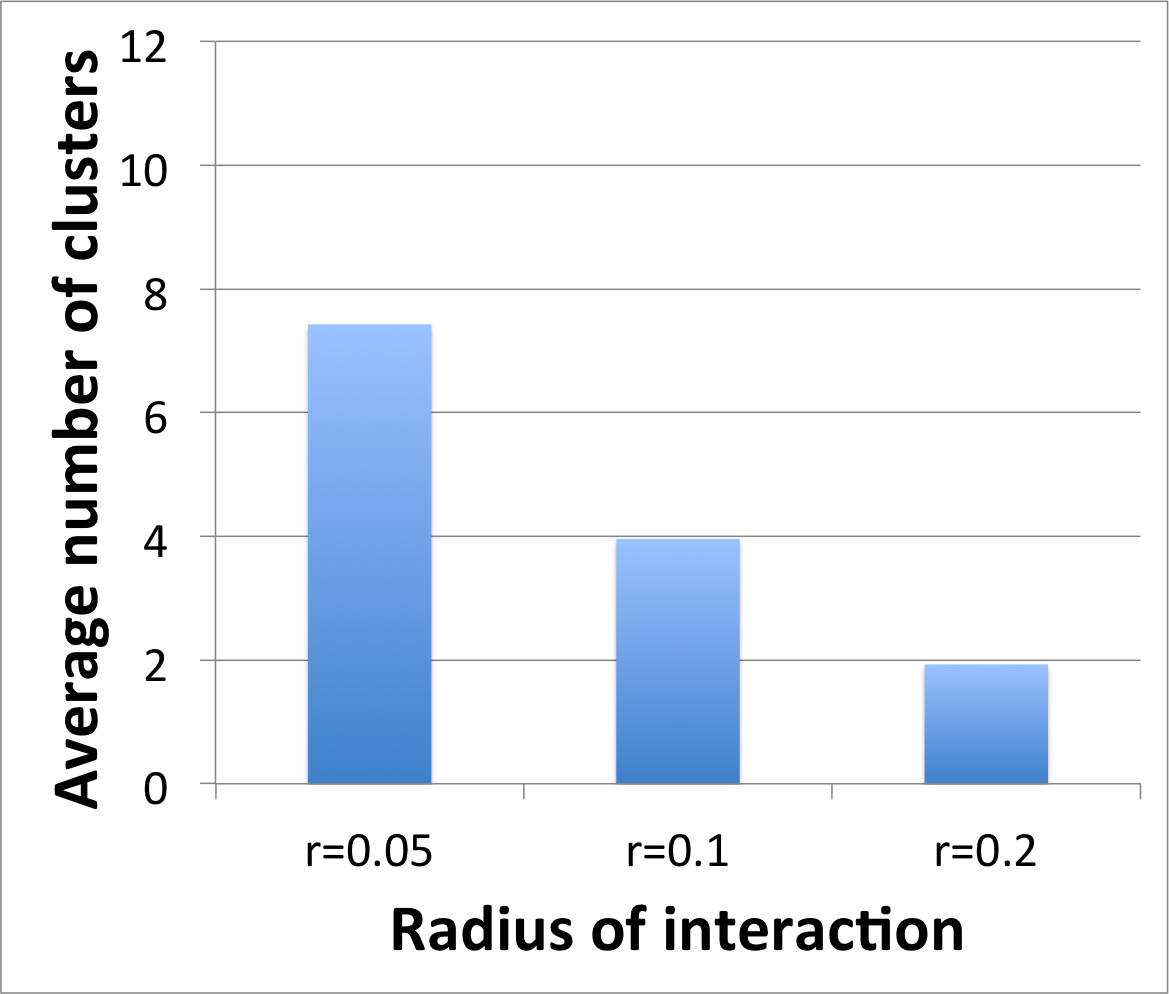}\label{fig:metricconnection}}\hskip 0.3in 
           \subfloat[Topological]{\includegraphics[width = 0.26 \textwidth, height=0.24 \textwidth]{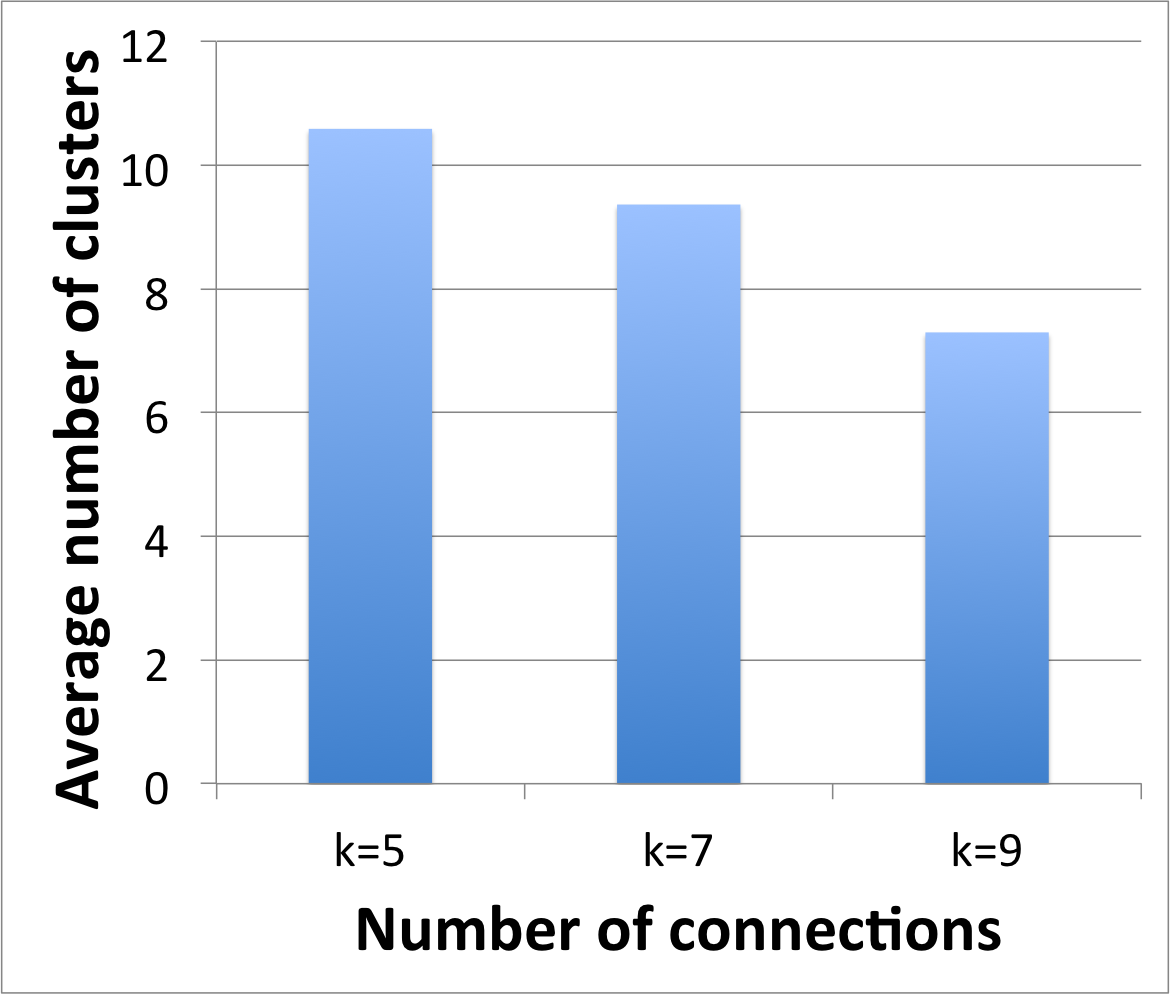}\label{fig:topoconnection}}\hskip 0.3in 
           \subfloat[Metric/topological]{\includegraphics[width = 0.30\textwidth, height=0.24 \textwidth]{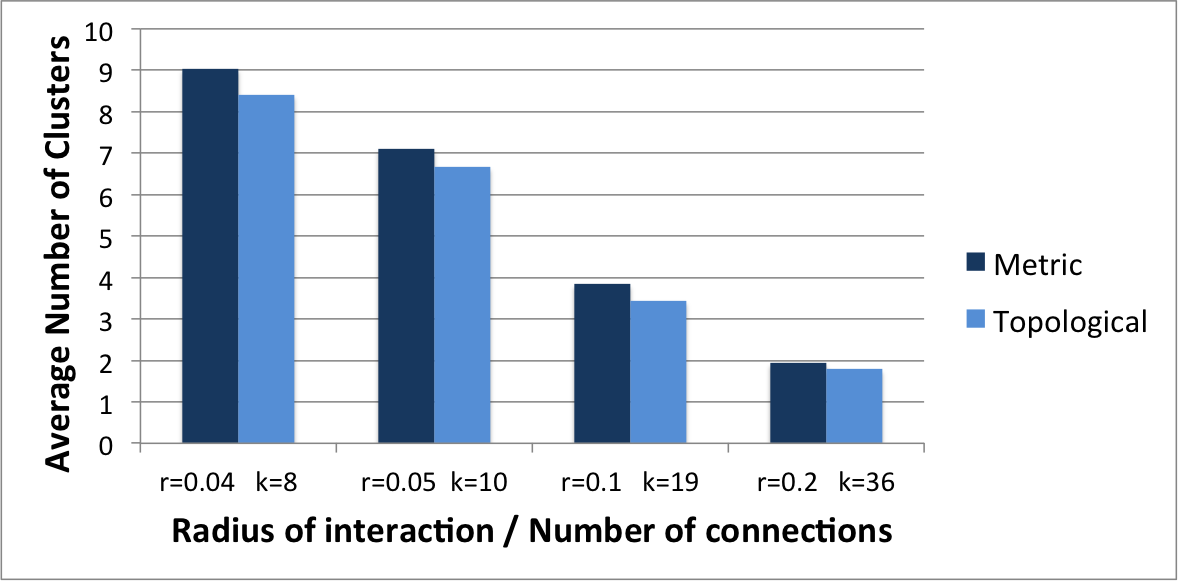}\label{fig:comparetopometric}}\hskip 0.3in 
        \caption{\small Average number of clusters of the asymptotic solution: (a) for different radii $r$ in the metric configuration, and (b) for different numbers of connections $k$ in the topological configuration. Each average was obtained over 100 simulations, in which 100 agents are initially distributed uniformly in the interval $[0,1]$. Figure (c) provides a comparison of the two networks, plotting side by side metric and topological configurations with the same initial average number of connections per agent.}
\label{fig:clustering}
\end{center}
\end{figure}

To compare the different interaction networks, we ran simulations for the well-established one-dimensional HK model, see Figure \ref{fig:clustering}.  
 Recent results \cite{BCC08} proposed the idea that topological interactions (with the 5-7 closest neighbors) is an effective way for birds to ensure group cohesion and to escape predators.
 Figures \ref{fig:metricconnection} and \ref{fig:topoconnection} show the average number of clusters of the asymptotic solution of the HK-model \eqref{eq:HK} respectively with metric interaction \eqref{eq:Nimetric} and with topological interaction \eqref{eq:Nitopo}, for a group of 100 agents. Notice that consensus is not reached for small radius of interaction ($r\leq 0.2$) or a small number of neighbors ($k<10$), but instead the group tends to cluster in several subgroups. As expected, the number of clusters decreases as the network connectivity increases.
% Contrary to the proposed ansatz, topological interaction with a small number of neighbors ($k<10$) is not very effective in ensuring consensus, which would imply group cohesion. 
 Figure \ref{fig:comparetopometric} shows that with the same initial number of connections, both interaction networks perform similarly.

\begin{figure}[h!]
\begin{center}
           \subfloat[$r=0.1$ ]{\includegraphics[ height=0.24 \textwidth, trim={10 10 50 40},clip=true]{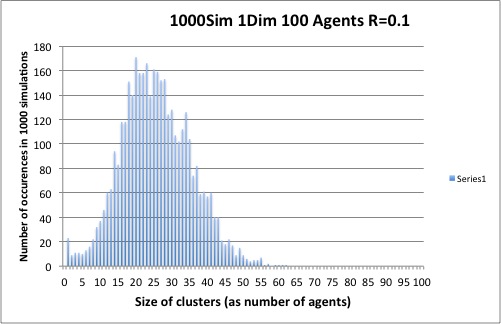}\label{fig:ClusterDistr01}}\hskip 0.3in 
           \subfloat[$r=0.2$]{\includegraphics[ height=0.24 \textwidth, trim={10 10 50 40},clip=true]{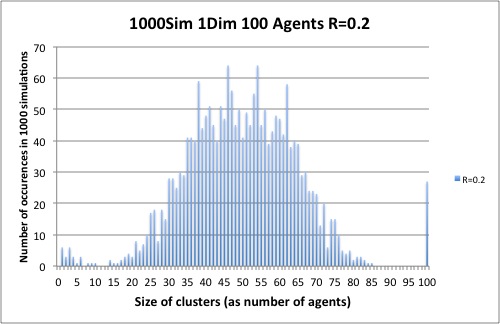}\label{fig:ClusterDistr02}}\hskip 0.3in 
\caption{\small Distribution of the asymptotic clusters' sizes in 1000 simulations of the one-dimensional HK model with 100 agents and metric interaction. Initially the agents are distributed uniformly in the interval $[0,1]$. Figure (a) was obtained with an interaction radius $r=0.1$ and Figure (b) with $r=0.2$. In the case $r=0.2$, consensus was reached in 28 simulations. Furthermore, the shape of the distribution suggests that some cluster sizes are more frequent than others.}
           \label{fig:ClusterDistr}
\end{center}
\end{figure}

\begin{figure}[h!]
\begin{center}
\includegraphics[ height=0.3 \textwidth, trim={10 5 10 40},clip=true]{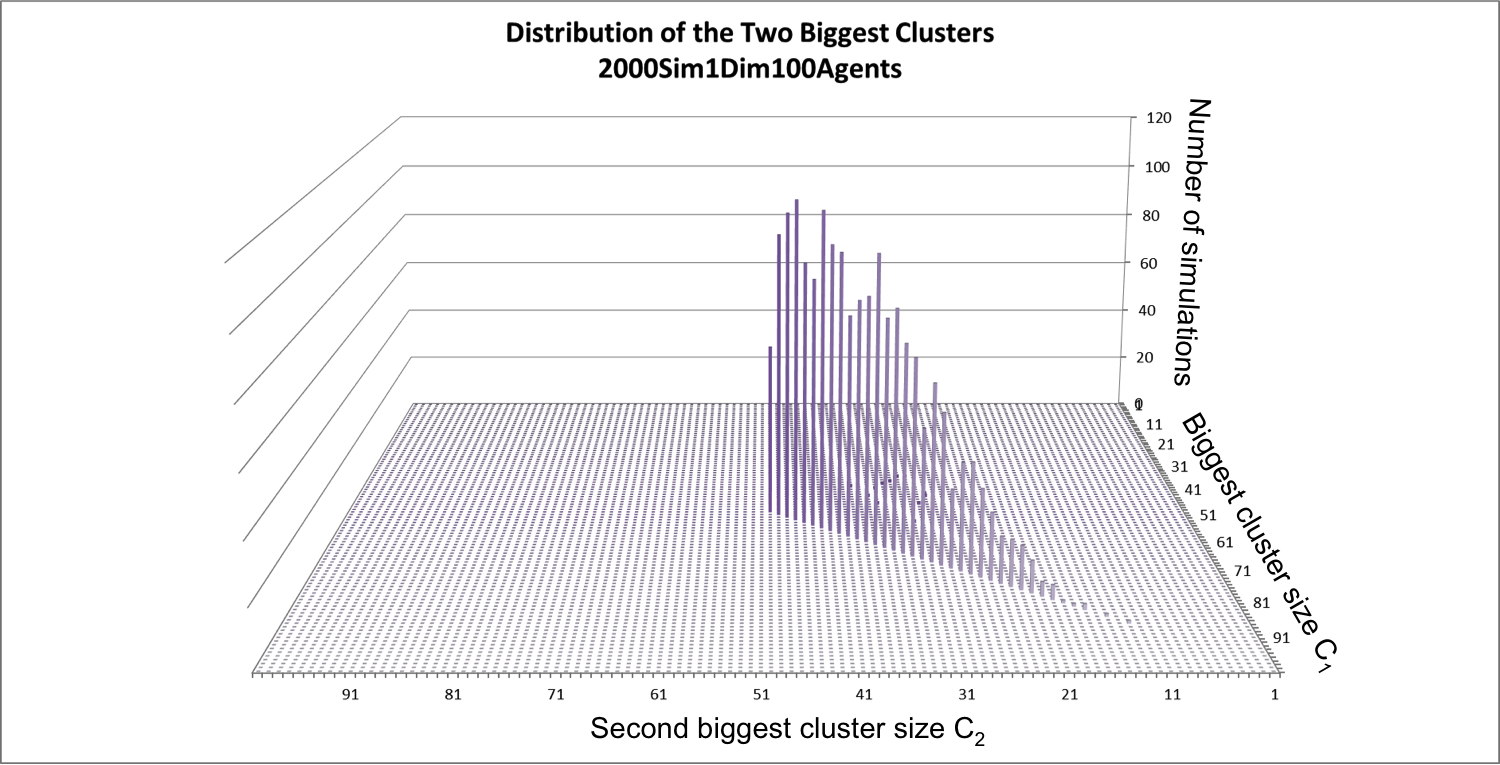}
\caption{\small Distribution of the two biggest asymptotic clusters' sizes in 2000 simulations of the one-dimensional HK model with 100 agents and metric interaction ($r=0.2$). Initially the agents are distributed uniformly in the interval $[0,1]$. The conformation $(C_1,C_2)=(50,50)$ is obtained in 60 cases, whereas the conformation $(C_1,C_2)=(51,49)$ is obtained in 100 cases. The most likely conformation is $(C_1,C_2)=(53, 47)$, obtained in 113 cases. There is a low likelihood of having $C_1-C_2>20$.}
           \label{fig:3dClusterDistr}
\end{center}
\end{figure}

In order to illustrate the differently stable equilibrium conformations, we ran simulations with the one-dimensional HK system, plotting the distribution of the asymptotic clusters' sizes (see Figure \ref{fig:ClusterDistr}). We observed that some conformations are statistically more frequent than others. For instance, in 1000 simulations of the HK dynamics of a group of 100 agents with metric interaction and an interaction radius $r=0.2$, clusters of 38, 46, 54 and 62 agents are the most frequently obtained (Fig.\ref{fig:ClusterDistr02}). 
Notice that if $r=0.2$ and the agents are distributed in the interval $[0, 1]$, there can be at most 4 clusters. We show that in the conditions of the simulations of Fig. \ref{fig:ClusterDistr02}, in most cases the agents are asymptotically distributed in 2 clusters. Figure \ref{fig:3dClusterDistr} shows the size distribution of the two biggest clusters $(C_1,C_2)$ over 2000 simulations. The peaks are mostly distributed along the line $C_1+C_2=100$, which means that in most simulations an equilibrium of 2 clusters is reached. Observe that it is less likely to reach an exactly equal distribution of agents between those two clusters than it is to have a slightly unbalanced distribution. The probability of having a very unbalanced distribution decreases with the imbalance.  

\begin{figure}[h!]
\begin{center}
           \subfloat[Evolution of the positions]{\includegraphics[height=0.25 \textwidth, trim={0 0 0 25},clip=true]{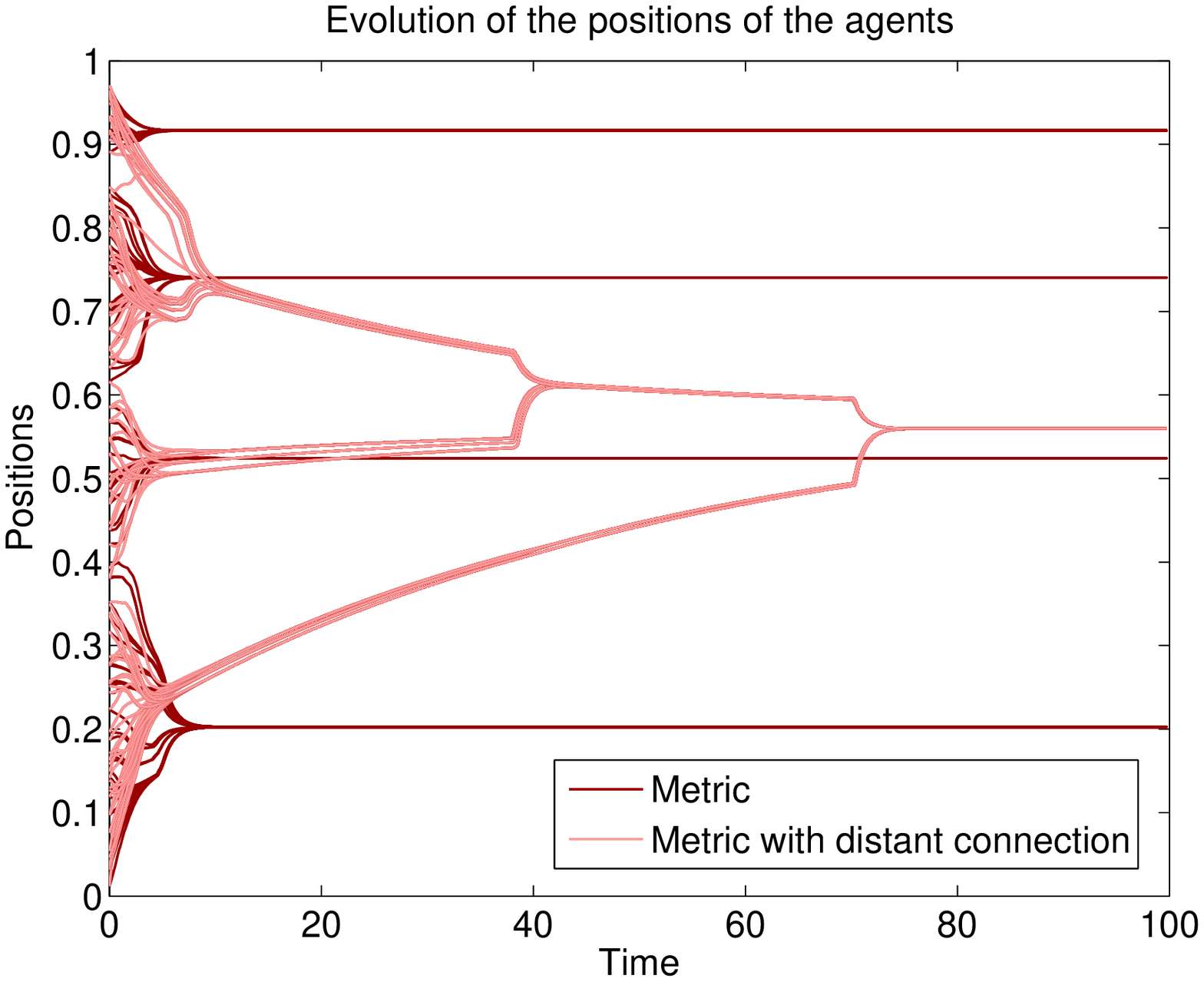}\label{fig:Evcousinpos}}\hskip 0.1in 
           \subfloat[Evolution of the number of edges]{\includegraphics[height=0.25 \textwidth, trim={0 0 0 25},clip=true]{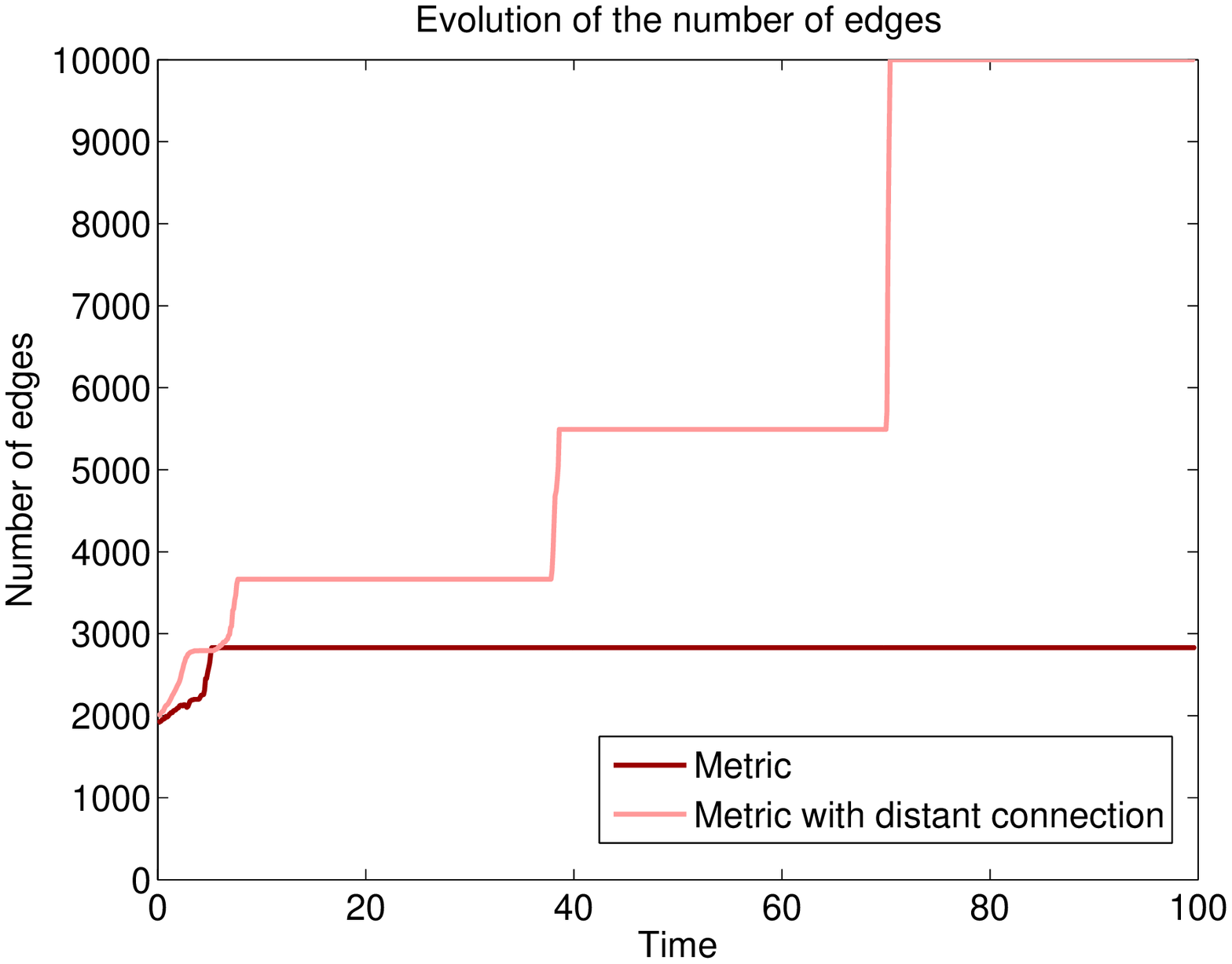}\label{fig:Evcousinedge}} 
        \caption{\small Effect of distant connections in convergence to consensus in the HK model with $r=0.1$. Figure (a) shows the evolution of positions in the metric case with only local interactions or with one added distant connection chosen uniformly (i.e. $a=0$), resulting respectively in clustering or consensus. Figure (b) shows the evolution of the number of edges.}
\label{fig:Evcousin}
\end{center}
\end{figure}

\noindent{Long-range connection}
We verified the effectiveness of long-distance connections in enhancing consensus for social dynamics. For each agent, a distant connection selected uniformly among the other agents was added to each agent's local interactions (see Figure \ref{fig:Nilong}). Added to metric interactions, the distant connection almost always lead to consensus. Figure \ref{fig:Evcousin} shows the improved convergence to consensus when adding an additional distant connection in the HK model. Figure \ref{fig:Evcousinpos} shows the evolution of positions with and without an added distant connection. Figure \ref{fig:Evcousinedge} shows the evolution of the total number of edges of the network. When distant connections are added, the system asymptotically reaches consensus, and the graph becomes fully connected, i.e. $\EE = \VV\times\VV$ so that $\mathbf{card}(\EE) = N^2$.

\begin{figure}[h!]
\begin{center}
           %\subfloat[Evolution of positions]{\includegraphics[width = 0.32 \textwidth, height=0.25 \textwidth]{Krause_metric_comparison.png}}\hskip 0.1in
           \subfloat[Time to consensus (metric)]{\includegraphics[width = 0.32 \textwidth, height=0.25 \textwidth]{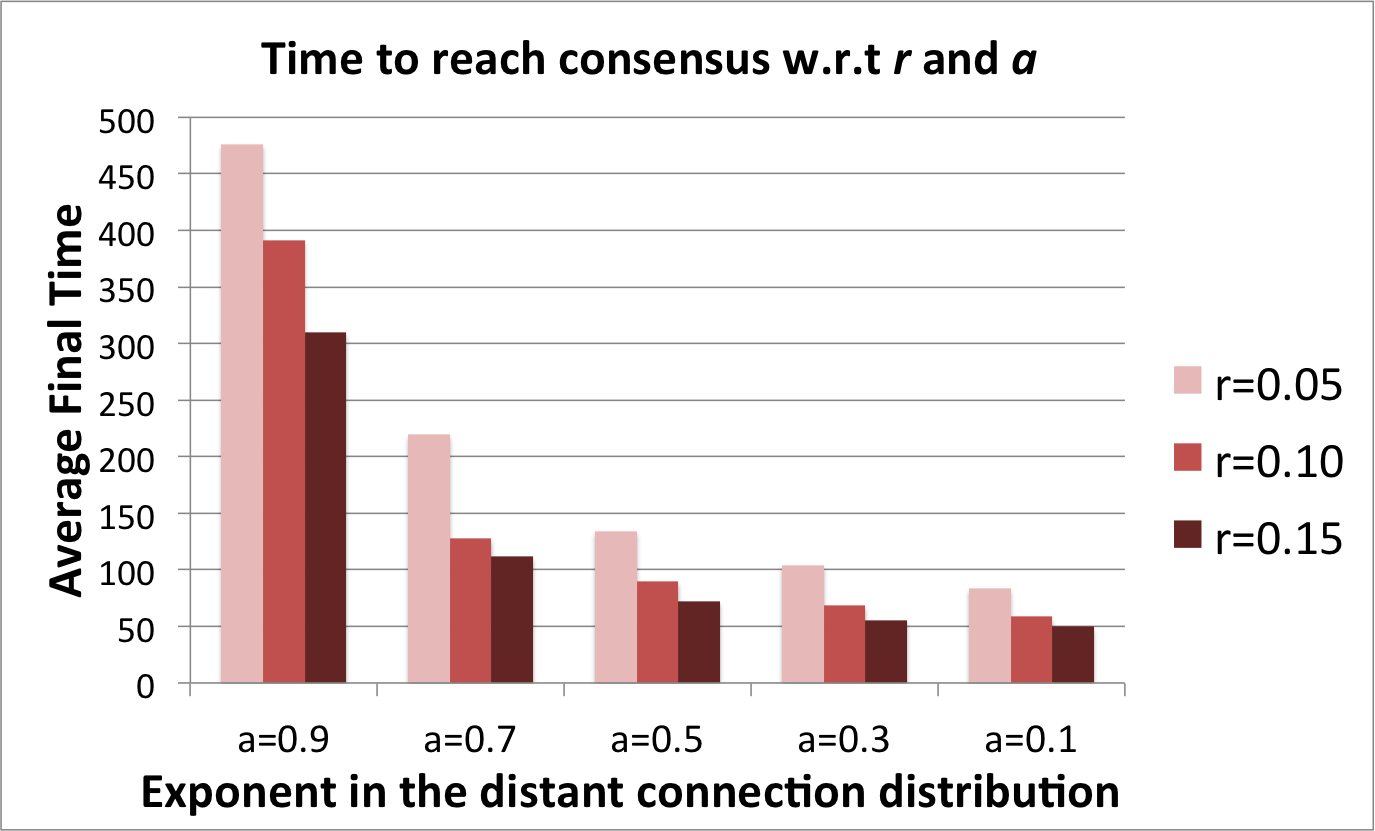}\label{fig:cousinmetric}} \hskip 0.1in 
           \subfloat[Clustering (topological)]{\includegraphics[width = 0.32 \textwidth, height=0.25 \textwidth]{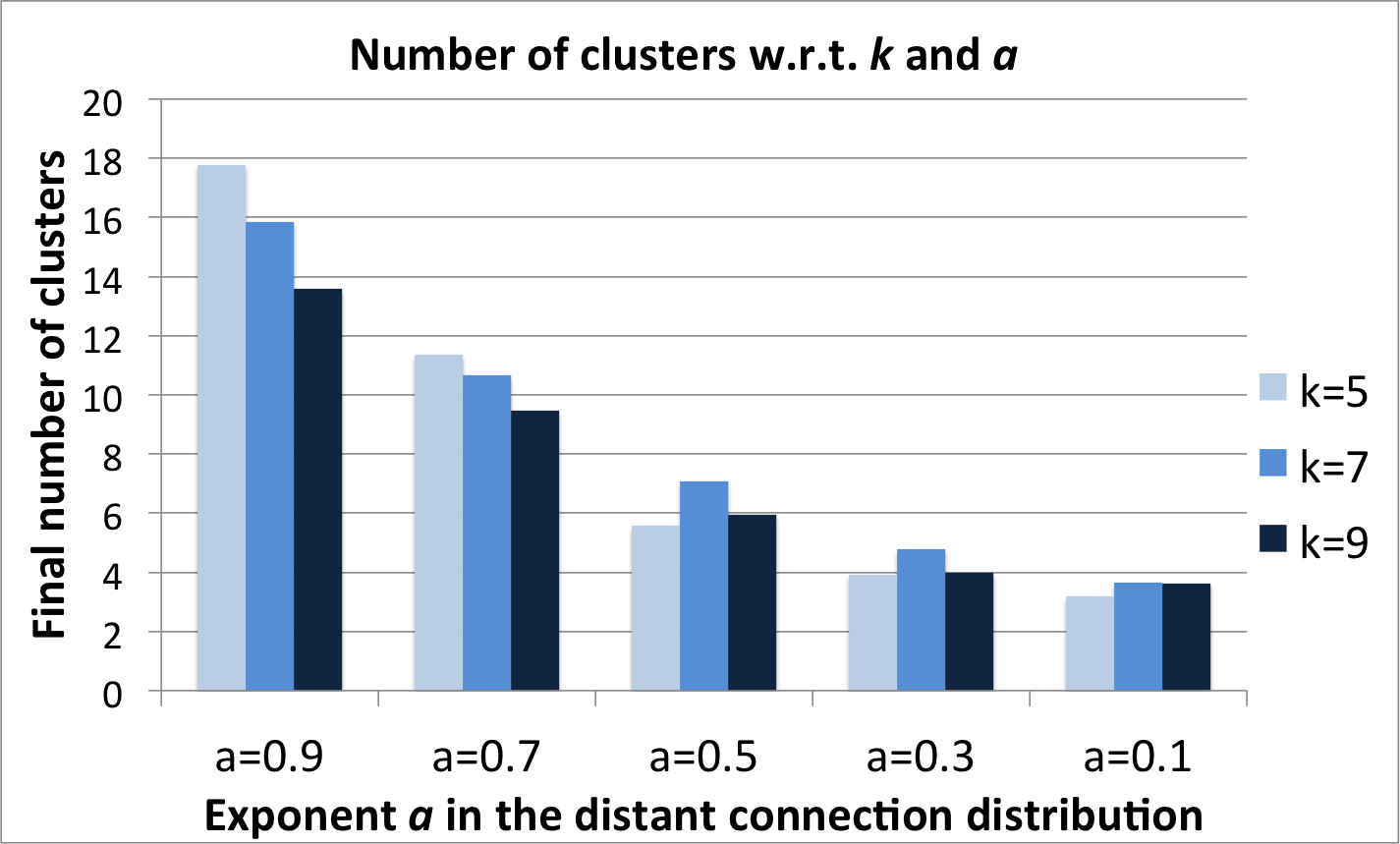}\label{fig:cousintopo}} 
        \caption{\small Effect of distant connections in convergence to consensus in the HK model. 
        %(a) Evolution of positions in the metric case resulting in clustering (with only local interactions, thin blue lines) or consensus (with one added distant connection, green thick lines) ; 
        Figure (a) shows the decrease of the time necessary to reach consensus by adding a distant connection (metric case). Since consensus is reached only asymptotically, time to consensus was defined as the time necessary for all agents to be within a sphere of given radius $\epsilon$. Figure (b) shows the decrease of the final number of clusters by adding a distant connection (topological case).}

\end{center}
\end{figure}

We then studied the effect of the probability with which the distant connection is chosen among all the graph edges. More specifically, we penalize the increase in distance between agents by choosing the distant connection with a probability proportional to 
 $\rho^{-a}$, where $a\in (0,1)$ and $\rho$ is the distance between agents.
With local metric interaction, adding such a distant connection almost always leads to consensus. With topological interaction, consensus is not always reached but the number of final clusters is significantly reduced.
The more biased the choice of distant connection is towards distant neighbors (i.e. the smaller the parameter $a$), the faster consensus is achieved in the metric case (Fig. \ref{fig:cousinmetric}) or the fewer clusters are obtained in the topological case (Fig. \ref{fig:cousintopo}). 

%Other than differences between metric and topological interactions, we will discuss the following aspects of interaction networks:
%\begin{itemize}
%\item Scale-free networks (such as the internet, social networks, or even cells) show a great degree of robustness by having a small number of high degree nodes \cite{AJB99, AJB00}. This ability to communicate despite a high rate of errors is counter-balanced by the network's vulnerability to attacks on these central nodes. 
%\item Kleinberg worked on the \textit{small world} problem, showing that adding long-distance connections to local ones greatly reduces the network's diameter and facilitates the spread of information \cite{K00}. 
%\end{itemize}

%% file: Networks-InteractionPotential.tex
 In equations \eqref{eq:syst-gen} and \eqref{eq:syst2-gen}, the interaction coefficient $a_{ij}$ can be defined as a function of the distance between the agents $i$ and $j$: $a_{ij}=a(\|x_i - x_j\|)$.
  % A common choice is to have the interaction potential between two agents $x_i$ and $x_j$ depend on the distance between the agents $ $.  
  The interaction potential $a$ can be chosen to be \textit{homophilious} or \textit{heterophilious}, a terminology used by Motsch and Tadmor \cite{MT14}.  Both of these interaction potentials are functions of the distance between two agents.
  \begin{itemize}[nolistsep]
  \item $a(\cdot)$ is a \textit{homophilious} interaction potential if it is a decreasing function of the distance between agents.
  \item $a(\cdot)$ is a \textit{heterophilious} interaction  potential if it is an increasing function of the distance between agents.
  \end{itemize}
  
 A homophilious interaction potential is appropriate when aiming to model behavior that has a strong interaction between two agents which are close together, or two opinions that are similar.  If a bird uses sight to maintain proximity to neighboring birds (a flock), then cohesiveness between birds will depend on distance and visual acuity.  A study of non increasing interaction functions in the CS model is given in \cite{CS07, HT08}.  Homophilious interaction potential is intuitive in the sense that, for agents, organisms, or opinions to influence each other, they must not be too distant.
 
 A heterophilious interaction potential agrees with the phrase, ``opposites attract''.  When two agents are very different from each other, they tend to have strong influence on the other.  Motsch and Tadmor show a counter intuitive result \cite{MT14}: heterophilious interaction potential increases clustering behavior.  Particularly, for the long term behavior of the system, the number of clusters of agents will decrease as the heterophilious interactions strengthen.  Sufficiently strong heterophilious interactions will drive the number of clusters to one, which is a consensus.
 
A more detailed model implements multiple interactions, such as short range repulsion and long range attraction.  This describes behavior where agents avoid collision but otherwise converge.  In this kind of model, there is cohesion among agents, but as soon as they come too close to each other, they will move apart.  This behavior is present in animal groups \cite{KR02,S06}, schools of fish \cite{HW92,P02}, and is used to model human crowds of pedestrians \cite{CPT11}.

\textit{Anisotropic interactions} are those that depend on an orientation of an agent relative to other agents.  For example, an animal may mostly receive information from its field of vision.  In this case, the visual space of an agent must be considered.  An animal may easily recognize animals in front, as opposed to behind.  A study of how these anisotropic interactions affect the structure of animal groups can be found in \cite{CFP11}.

%% file: StateSpace-Opinion.tex
%{\bf Discrete states.}

%The Sznajd model, see \cite{SZ}, is based on the Ising model for ferromagnetism in statistical mechanics. In this model opinions are discrete variables $x_{i}$ taking value $\pm 1$. The interactions are governed by
% two basic rules: the ``ferromagnetic'' interaction (that is, if $x_{i} = x_{i+1}$ then at the next step adjacent agents will satisfy with a given probability $x_{i-1} = x_{i} = x_{i+1}= x_{i+2}$) and the ``antiferromagnetic'' interaction 
% (if $x_{i} = - x_{i+1}$ then an antisymmetric pattern forms $- x_{i-1} = x_{i} = - x_{i+1}= x_{i+2}$). The model has been extended to higher dimensional opinion and complex network topologies.  The motivation for this model comes form the postulate that ``agreement  generates agreement'', that is, if two agents reach a consensus then all agents directly connected to them are induced to agree. In other words, in Sznajd model, the opinion flows out from a group of agreeing agents (Social Validation). 
 
%The Sznajd model belongs to the class of binary-state opinion dynamics model. Another example 
% of opinion formation model in this class is the well known  Voter model that can be seen as a generalization of the Sznajd model (\cite{behera}). 

%{\bf Compact Manifolds.}
To describe the slow and continuous evolution of opinions, several models adapted consensus algorithm on Euclidean spaces, such as the HK model described previously (Section \ref{Sec:Overview-order1}). 
%This is the case of the well known Hegselmann--Krause bounded confidence model \cite{HK}, for instance, which is a model of herding of opinions in a $N$-agents system in which the position $x_{i}$ of the agent $i$, representing its opinion and taking values in an interval of $\mathbb{R}$, changes according to the distance from other agents $x_{j}, j \neq i$, rescaled by an interaction coefficient $a_{ij}$ accounting for the weight given to the opinion of agent $j$ by agent $i$. With these notations the opinion of agent $i$ evolves, in the continuous-time version of the model according to
%$$
%\dot x_{i} = \sum_{|x_{i}-x_{j}|<1} a_{ij} (x_{j} - x_{i}).
%$$ 
%The rationale for the bounded confidence is that it is unlikely for one agent to be influenced by another one whose opinion is too far from its own.  This kind of interaction gives rise to clusters of opinions (see for instance~\cite{BHT}). We also mention the bounded confidence model by Deffuant Model, see~\cite{deffuant} in which the opinions belong to real intervals too but the pairs of agents in interaction are chosen randomly.
The dynamics of consensus models on Euclidean spaces are, at least locally, linear. This may be a limitation in representing the complex behavior of opinions. Indeed the only equilibria of the system are clusters of consensus (see~\cite{BHT}). This is one of the main issues determining  a lack of connections with real life examples as pointed out by Sobkowicz \cite{sob}. 

Recently there has been a growing interest in designing consensus algorithm on nonlinear manifolds. The motivation comes from engineering applications, indeed oscillators evolve on the circle $S^{1}$, satellite altitudes evolve on the special orthogonal group $SO(3)$ and ground vehicles on the euclidean groups $SE(2)$ or $SE(3)$.
The first model in this direction is the Kuramuto model~\cite{kuramoto} on the sphere $\mathbb{S}^1$ which attracted a wide interest of researchers over the last 30 years,  motivated by its  connection with the problem of synchronizing a large population of harmonic oscillators - see the survey by Strogatz \cite{Kuramuto-survey}. Other possible applications
were studied by Hopfield \cite{hopfield1982} and Vicsek et al. \cite{VCABC95}. Lately, convergence analysis for adapted versions of the Kuramoto model on the circle have been thoroughly  studied in a series of papers by D\"orfler, Chertkov and Bullo \cite{DMB13}
Scardovi, Sarlette and Sepulchre \cite{SSS2007},
Sepulchre, Paley and Leonard \cite{SPL2007,SPL2008}.

 A first effort in studying consensus dynamics on more general manifolds  has been made  in~\cite{sarsep2009} by Sarlette and Sepulchre who looked at,
% studied opinion dynamics on a wider class of manifolds including, 
 among others, the special orthogonal group $SO(n)$, the Grassmann manifold, and $\mathbb{S}^1$ (see also \cite{sarlette2009} and \cite{sepulchre2011} for a survey on this topic).
Consensus problems on general manifolds present an inherent difficulty: in order to move towards a given point (for instance the weighted average of its neighbors' positions), an agents must follow the geodesics of the manifold, which are well defined only locally. Not only can geodesics not be unique on a global scale, but their computation can be extremely challenging. 
One way around this difficulty is to consider the embedding of the manifold $M$ into a Euclidean space $E$ (for instance $E = \R^d$). Using the embedding,
 these models are mainly based on the projection of linear consensus dynamics on the tangent space to the manifold $M$. Namely, given $N$ agents, their opinions $x_{i} \in M$ evolve according to:
\begin{equation}\label{eq:consensusonmanifolds}
\dot x_{i} = \Pi_{x_{i}}\left(\sum_{j=1}^{N} a_{ij} (\hat x_{j} - \hat x_{i})\right), \quad \mbox{ for } i=1,\ldots,N,
\end{equation}
where $\hat{x}_i$ denotes the embedding of $x_i$ in $E$ and $\Pi_{x}(y)$ is the projection of $y$ onto the tangent space to $M$ at $x$.
The dynamics for these systems inherit locally the structure of the linear case and convergence results rely mainly on 
consensus algorithms for linear systems (as, for example, the one by Tsitsiklis \cite{tsitsiklis1984}, Jadbabaie, Lin, and Morse \cite{jadbabaie2003}, Moreau \cite{moreau2004,moreau2005}, Blondel, Hendrickx, Olshevsky and Tsitsiklis
\cite{BHT}, Olfati-Saber and Murray \cite{olfati2007}, etc.)  but this is no longer possible  for global convergence analysis since  the considered manifolds are in general not globally convex.
As in linear algorithms consensus is an equilibrium of the system. In Euclidean 
spaces, if the interaction graph associated with the interaction coefficients $a_{ij}, i, j=1,\ldots,N$ is strongly connected, then
the system always tends to consensus. In nonlinear manifolds, consensus configurations become graph-dependent.

Systems on compact manifolds show more diverse kinds of equilibrium configuration, for instance the anti-consensus, in which 
 each state is furthest from the mean of its neighbors, so that as a result the opinion spread over the entire manifold. This phenomenon is sometimes called \emph{balancing}, in opposition with the term \emph{synchronization} used to describe consensus on the circle \cite{SSS2007}.

Recently Caponigro, Lai and Piccoli proposed a nonlinear opinion formation model on the $d$-dimensional sphere $\mathbb{S}^{d}$~\cite{CLP15}.
The rationale for the sphere $\mathbb{S}^{d}$ is that, as mentioned, opinions are subjected to a quantization phenomenon when measured. We can imagine that, at the instant of measurement (elections, polls, interviews, etc.), opinions take only two values (yes/no, left/right, Democratic/Republican, liberal/conservative, for/against, etc.), so that every component of the vector 
$x_{i} = (x_{i}^{(1)}, \ldots, x_{i}^{(d+1)})$ takes a positive or a negative value. In particular $x_{i}$ belongs to $\sqrt{d+1}\;\mathbb{S}^{d}$.  The manifold $\mathbb{S}^{d}$ is a  mathematical abstraction to describe the dynamical evolution of the opinions on a continuous (i.e. non-discrete) set. Moreover opinions on different topics are usually interconnected:
economic policy attitudes and candidate choice in political elections; opinion formation and economical condition; opinion on research funding and religious or ideological beliefs. System~\eqref{eq:consensusonmanifolds} on the sphere shows new  kind configurations with respect to the one observed on Lie Groups. Beside consensus also antipodal and polygonal equilibria appear in the model (that can be seen as balanced configuration). Furthermore a  configuration typical of this system, called \textit{dancing equilibrium}  is shown. In this configuration the mutual distances between opinions are in equilibrium but the system may evolve.

%% file: StateSpace-Animal.tex
%We can study second-order models with more general manifolds as state-space, such as what was done for opinion-formation models in \cite{CLP15}.
%See also \cite{DM08} and \cite{SBS10}.

Standard second-order social dynamics systems such as the CS dynamics \eqref{eq:CS} evolve in the Euclidean space $\R^{2Nd}$ where $N$ is the number of agents and $2d$ the dimension of the state space for the position and velocity (typically $d=2$ or $d=3$). However, similarly to opinion dynamics models, some applications require more complicated state-spaces. 

One of the main difficulties in modeling opinion formations is the lack of reliable methods to measure opinions. A classical problem in sociology is to design interviews not affecting opinions, i.e. questions not influencing answers. Purely open questions do not exist and, moreover, it is very hard to collect data from open answers. On the other hand, closed questions induce quantization on the answers: opinions collapse on discrete sets representing the possible answers to a closed question.  
This is the rationale to design models in which the  initial and final opinions, in an opinion formation process, take value in a discrete set as in well known Sznajd model and Voter Model, for instance.
As seen in section \ref{Sec:Overview-order2}, the Sznajd model belongs to the class of binary-state opinion dynamics model. It is based on the Ising model for ferromagnetism in statistical mechanics \cite{SZ}. In this model opinions are discrete variables $x_{i}$ taking value in the space $\{-1,1\}$.  It has been established that there are two possible equilibria for this model: \textit{ferromagnetism}, in which all agents have the same spin, and \textit{anti-ferromagnetism}, in which agents have alternate spins.

%The Sznajd model belongs to the class of binary-state opinion dynamics model. Another example 
% of opinion formation model in this class is the well known  Voter model that can be seen as a generalization of the Sznajd model (\cite{behera}). 
Another classic example is the Vicsek model \cite{VCABC95} in which every particle's velocity is assumed to have constant norm, so that each particle is represented by its two-dimensional position and the angle of its velocity. This model allows to study clustering and orientational order, two patterns commonly observed in various biological systems such as animal groups or bacteria.
As a variation upon the Vicsek model, Motsch and Degond designed the persistent turning walker model in order to study fish motion, where the velocity is also assumed to have constant norm $c$  \cite{DM08}. The variables are the two-dimensional position of the fish's centroid $x\in\R^2$, the velocity angle $\theta\in\R /2\pi\mathbb{Z}$ and the curvature of the trajectory $\kappa\in\R$. The trajectories are described by the stochastic differential equations:
\begin{equation*}
\begin{cases}
\dot{x} = c\tau(\theta) \\
\dot{\theta} = c\kappa \\
d\kappa = -a\kappa dt+b dB_t
\end{cases}
\end{equation*}
where $\tau(\theta)=(\cos\theta,\sin\theta)$ is the direction of the velocity vector, $dB_t$ is the standard Brownian motion, $a$ is a relaxation frequency and $b$ quantifies the intensity of the random curvature jumps.
 The dynamics of the curvature of the trajectory reflect the antagonistic effects of its tendency to relax to a straight line and of the random jumps observed in fish behavior.

%% file: Control.tex
% Sparse control, stabilization, optimal control, emergence of consensus

Many works have explored ways of controling social dynamics systems. Control can be of great use, for instance in applications to robotics, for rendez-vous problems. One can aim to control the system to: reach consensus in the state space or alignment in the velocity space \cite{CFPT13, CFPT15};
reach a predetermined desired position or velocity \cite{PPS15};
keep the agents as far from each other as possible, to avoid Black Swan type phenomena where consensus can lead to market collapse (in applications to economics). 

Control is a particularly challenging problem due to the high dimensionality of the systems.
One can either control the high-dimensional discrete system \cite{PPS15}; 
resort to mean-field control \cite{HPS15,FPR14};
act on the network to exploit its intrinsic properties (such as symmetry) \cite{RJME09}.
 Control often leads to separating the group into a set of controlled leaders and a set of uncontrolled followers.

% \bigskip
% \bigskip
% {\color{red}\bf
% 
% ------------------------------------------------------------------------
% 
% \medskip
% 
% NB: I have made a (very) draft redaction, because for the moment the rest of the file is empty. Smoothing and proper referencing will be strongly required. 
% 
% There are tons of possible references: hereafter I have only given names of authors. You will please choose what to cite.
% 
% ------------------------------------------------------------------------
% 
% }

%Objectives:\\
%Model \textit{self-organization} and \textit{consensus emergence} in a group of agents evolving in time according to mutual interactions. There are:
%\begin{itemize}
%\item models in finite dimension (Cucker-Smale)
%\item models in infinite dimension (mean-field limit)
%\end{itemize}
%We want to promote \textit{organization via intervention} (social forces), enforce or facilitate \textit{pattern formation} or \textit{convergence to consensus}:
%\begin{itemize}
%\item add a control in the model
%\item design an appropriate control strategy
%\end{itemize}

\subsection{Finite dimension} \label{Sec:ControlFinite}

We start by presenting various control techniques related to finite-dimensional models. 

\subsubsection{Consensus protocols}
Since the first 2000s consensus in multiagents systems has been seen also as a distributed control problem (see for instance Jadbabaie, Lin and Morse \cite{jadbabaie2003}, Olfati-Saber, Fax and Murray \cite{olfati2007} and also Tsitsiklis \cite{tsitsiklis1984}).
The problem in this framework is to find a feedback control, or \emph{consensus protocol},
$u(x)=  (u_1(x),\ldots,u_N(x))$ assigning dynamics to the $i$-th agent
$$
\dot x_i = u_i(x)
$$
for $i=1,\ldots,N$ such that every solution tends to a consensus configuration $x_1= \dots = x_N$.
The first results in this direction deal with linear consensus algorithm of the form
$$
u_i(x) =\sum_{j=1}^N a_{ij} (x_j-x_i),
$$
and show sufficient conditions guaranteeing asymptotic consensus under minimal
connectivity assumptions on the communication graph associated with the interaction coefficients $a_{ij}$, $i,j =1,\ldots,N$ (see for instance Moreau ~\cite{moreau2004,moreau2005}).

%\medskip
\subsubsection{Non-consensus}
A great interest has been given to the modeling of emergent behavior in animal groups and social dynamics (See 
Section \ref{Sec:Overview}). However self-organization is not always sufficient to ensure consensus of positions or alignment of velocities.
%At the same time a great interest has been given to the modeling of emergent behavior in animal groups and social dynamics (See 
%Section \ref{Sec:Overview}).  
%In order to describe the emergence of behaviors with different natures or complex phenomena, models designing focuses on multiagents systems in which a steady configuration, as consensus, is not reached asymptotically by any trajectory or, alternatively, in which several kind of equilibria are present. 
%
%A classical example is given by the well-known Hegselmann--Krause bounded confidence model~\cite{HK} in which the asymptotic behavior shows clusters of opinions. 
%, representing the  consensus for this second order model.
The following example shows initial conditions for which the Cucker-Smale system does not tend to alignment.

\begin{remark}
Consider the CS system ~\eqref{eq:CS} in the case of two agents moving in $\R$ with position and velocity at time $t$, $(x_{1}(t),v_{1}(t))$ and $(x_{2}(t),v_{2}(t))$. Assume that $a(x) = 2/(1+x^{2})$. Let $x(t) = x_{1}(t) - x_{2} (t)$ be the relative main state and $v(t) = v_{1} (t) - v_{2} (t)$ be the relative flocking parameter. Then~\eqref{eq:CS} reads
$$
\left\{
\begin{split}
\dot x &= v\\
\dot v &= -\frac{v}{1+x^{2}}
\end{split}
\right.
$$
with initial conditions $x(0)=x_{0}$ and $v(0)=v_{0} > 0$. 
The solution of this system can be found by direct integration, {as from $\dot v = - \dot x / (1+x^{2})$ we have}
$$
v(t) -  v_{0} = -\arctan{x(t)} + \arctan{x_{0}}.
$$
%If the initial conditions satisfy $|\arctan{x_{0}} + v_{0}|  < \pi/2$ then it is easy to see that the relative main state $|x(t)|$ is bounded uniformly by $\tan{(|\arctan{x_{0}} + v_{0} |)}$, otherwise we would have $v(t^*) =0$ for a finite $t^*$. The boundedness of $x(t)$ fulfills the sufficient condition on the states in  Proposition~\ref{propHaHaKim} for {consensus}.  If $ |\arctan{x_{0}} + v_{0}|  = \pi/2$  then the system tends to consensus as well, since $v(t) = \pm \pi/2 - \arctan{x(t)}$, depending on whether $\pm v_0 >0$ respectively: if $x(t)$ were unbounded then $\lim_{t \to \infty} x(t)= \pm \infty$, respectively, and necessarily we converged to consensus. If $x(t)$ were bounded then again by Proposition~\ref{propHaHaKim} we would converge to consensus. \\
%On the other hand, 
Whenever  the initial conditions satisfy  $|\arctan{x_{0}} + v_{0}| > \pi/2 $,  which implies $|\arctan{x_{0}} + v_{0}| \geq  \pi/2 +\varepsilon$ for some $\varepsilon>0$, the  flocking parameter $v(t)$ remains bounded away from $0$ for every time, since
$$
|v(t)|= |-\arctan{x(t)} + \arctan{x_{0}} + v_{0}| \geq | -\arctan{x(t)} + \pi/2 +\varepsilon| > \varepsilon, 
$$
for every $t>0$.
In other words, the system does not tend to flocking.
\end{remark}

\subsubsection{External control}
%Therefore it is natural to study whether is possible to enforce an emergent behavior when this does not arise from the free evolution. 
%
When flocking is not achieved by self-organization, it is natural to wonder whether it is possible to control the group to flocking by means of an external action. We are therefore concerned with \textit{organization via intervention}.
Since flocking is a steady configuration of the system, enforcing self-organization can be seen as an asymptotic stabilization
problem, which is classical in control theory and usually relies on Lyapunov design 
(see for instance Isidori \cite{IsidoriBook} or Sontag \cite{SontagBook}). 
Using these classical techniques it is easy to design a stabilizing feedback. To better understand the problem let us consider the 
Controlled CS model, introduced by Caponigro, Fornasier, Piccoli and Tr\'elat \cite{CFPT13, CFPT15}. Consider the control system
\begin{equation}\label{sys_CS_control}
\left\{
\begin{split}
\dot{x}_i(t)&=v_i(t) \\
\dot{v}_i(t)&=\frac{1} {N}\sum_{j=1}^N a(\Vert x_{j}(t) - x_{i}(t)\Vert )(v_j(t)-v_i(t)) + u_i(t)
\end{split}
\right.
\qquad\qquad i=1,\ldots,N
\end{equation}
with the bound on the control
\begin{equation}\label{eq:controlbound}
\sum_{i=1}^{N}\Vert u_{i}(t)\Vert  \leq M
\end{equation}
for a given $M>0$. 
Any $v\in\R^d$ can be written as $v=(\bar v,\ldots,\bar v)+v_\perp$.
In \cite{CFPT15} the authors proved that 
the feedback control defined by 
\begin{equation}\label{eq:CScontrol}
u(t) = - \alpha v_{\perp}(t),
\end{equation}
for $\alpha>0$, stabilizes the system to flocking (in infinite time) while satisfying condition \eqref{eq:controlbound} if $\alpha$ is small enough.
Indeed, $\dot{V} \leq -\frac{2}{N} \sum_i \langle {v_\perp}_i,u_i\rangle = -2\alpha V$.

\begin{remark}
The control \eqref{eq:CScontrol} acts on a large number of agents simultaneously. This is inconvenient for practical purposes, since it requires intensive instantaneous communications between all agents. In what follows we look at more economical controls that are active on as few components as possible at any instant of time. This leads to the concept of \textit{sparse control}.
\end{remark}

\subsubsection{Sparse stabilization}\label{Sec:Control-Finite-Sparse}

%Of course in practice it is desirable to have as few dogs as possible, meaning that we would like the control to have only few components that are active, at any instant of time. This leads to the concept of \textit{sparse control}.

The objectives of sparse stabilization are:
\begin{itemize}[nolistsep]
\item To design a \textit{sparse feedback control} steering ``optimally'' the system to flocking, with:
\begin{itemize}[nolistsep]
\item[(i)] a minimal amount of components active at each time: concept of \textit{componentwise sparse control}.
\item[(ii)] a minimal amount of switchings in time: concept of \textit{time sparse control}
\end{itemize}
\item To control the system to any prescribed flocking.
\end{itemize}

Our idea to promote sparsity is to use $\ell^1$ minimization, as in image analysis where it has become very popular.

Note that the (far from being sparse) feedback stabilizing control \eqref{eq:CScontrol} 
is solution of the minimization problem
$$
\min_{\sum_{i=1}^{N}\Vert u_{i}\Vert  \leq M} \left( \frac{1}{2N^2}\sum_{i,j=1}^N\langle v_i-v_j,u_i-u_j\rangle  \right)
=
\min_{\sum_{i=1}^{N}\Vert u_{i}\Vert  \leq M} \left( \frac{1}{N}\sum_{i=1}^N\langle {v_\perp}_i,{u_\perp}_i\rangle \right).
$$
Instead, we now consider the slightly modified minimization problem
$$
\min_{\sum_{i=1}^{N}\Vert u_{i}\Vert  \leq M} \left( \frac{1}{2N^2}\sum_{i,j=1}^N\langle v_i-v_j,u_i-u_j\rangle + \gamma(X) \frac{1}{N} \sum_{i=1}^{N}\Vert u_{i}\Vert   \right) ,
$$
where 
$$
\gamma(X) = \int_{\sqrt{X}}^{+\infty} a (\sqrt{2N}r)\, dr.
$$
Here, the use of the $\ell^1$ norm is to enforce sparsity, and the weight $\gamma(X)$ is used as a threshold implying that the control will switch off when entering the flocking region.
The optimal solution of this minimization problem is the \textit{componentwise sparse feedback control} $u^\circ$ defined as
\begin{itemize}[nolistsep]
\item if $\displaystyle\max_{1\leq i\leq N} \Vert v_{\perp_{i}}(t)\Vert  \leq  \gamma(X(t))^{2}$, then $u^\circ(t)=0$
\item if 
$\displaystyle \Vert v_{\perp_{j}}(t)\Vert  = \max_{1\leq i\leq N} \Vert v_{\perp_{i}}(t)\Vert  >  \gamma(X(t))^{2}$ (with $j$ be the smallest index) then 
$$
u^\circ_{j}(t) = - M \frac{v_{\perp_{j}}(t)}{\Vert v_{\perp_{j}}(t)\Vert },\quad\mbox{ and }\quad u^\circ_{i}(t) = 0 \quad \mbox{for every } i\neq j.
$$
\end{itemize}

\begin{figure}[h]
%\centerline{\includegraphics[width=7cm]{controltoconsensus.pdf}}
\begin{center}
\begin{tikzpicture}
\draw[scale = 0.5, thick,->] (0,0) -- (11,0) node[right] {$X_0$};
\draw[scale = 0.5, thick,->] (0,-0.1) -- (0,11) node[left] {$V_0$};
%\draw[scale = 1,domain=0:1,smooth,variable=\x,blue] plot ({\x},{ 1/(200)*
%((3.14)/2- atan(sqrt(\x)))^2 });
\draw[scale = 0.5,domain=0.01:10,smooth,variable=\x,black] plot[id=x,samples=200] ({\x},{4*((3.14)/2- rad(atan(sqrt(\x))))^2});
\fill [scale = 0.5, gray, domain=0:10, variable=\x,samples=200]
      (0, 0)
      -- plot ({\x}, {4*((3.14)/2- rad(atan(sqrt(\x))))^2})
      -- (10, 0)
      -- cycle;
\draw	(0,-0.1) node[above right] {\color{white} \emph{Flocking region}};
\draw	(0,0) node[below left] {0};
\draw (4,3) node[align = right, right] (start) {\textbullet $(x_0,v_0)$};
\draw (0.5,0.7) node[align = right, right] (end) {\textbullet $(x(t),v(t))$};
\draw[<-] (0.9,0.9) .. controls (1.4,1.5) .. (3.9,2.9);
\end{tikzpicture}
\caption{\small Control to flocking. The flocking region $\sqrt{V_0}\leq\frac{1}{\sqrt{2N}}(\frac{\pi}{2}-\arctan(\sqrt{X_0}))$ corresponds to the CS model \eqref{eq:CS} with parameter $\beta=1$, see Proposition \ref{propHaHaKim}.}\label{fig_cont_CS}
\end{center}
\end{figure}
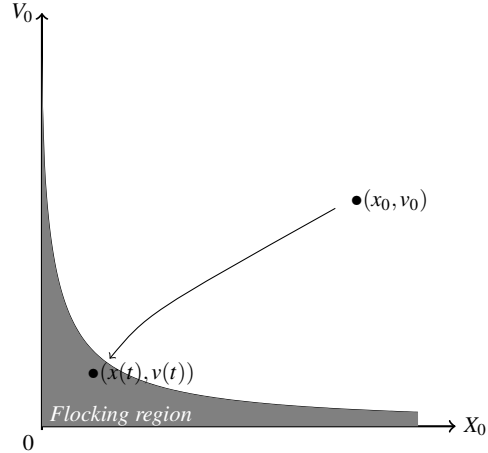

\begin{theorem} \cite{CFPT13, CFPT15}
The sparse feedback control $u^\circ$ stabilizes the system to flocking.
\end{theorem}
Indeed, we have
$$
\dot V \leq \frac{2}{N}\sum_i\langle {v_\perp}_i,u^\circ_i\rangle
= -2\frac{M}{N} \Vert {v_\perp}_j\Vert
$$
with $\Vert v_{\perp_{j}}\Vert  = \max_{1\leq i\leq N} \Vert v_{\perp_{i}}\Vert  \geq \sqrt{V}$ implying that $\dot V\leq -2\frac{M}{N}\sqrt{V}$, 
hence any trajectory enters in finite time the flocking region, and then we take $u=0$ (forever), as illustrated on Figure \ref{fig_cont_CS}, thus letting the trajectory naturally converge to flocking.
Note that, alternatively, one can choose not to switch off the control (even when one has entered the flocking region): in that case, the trajectory reaches flocking within finite time, because $\sqrt{V(t)}\leq \sqrt{V(0)}-2\frac{M}{N}t$.

\begin{remark}
By construction, this feedback control is componentwise sparse. However, it is not necessarily time sparse: it may \textit{chatter}. Indeed, the strategy described above consists of focusing on the agent that is the farthest possible from the mean, in order to steer it closer to the mean. But that agent may change as time evolves and oscillations may appear.
In order to avoid possible chattering in time, \cite{CFPT13, CFPT15} have implemented the classical sample-and-hold procedure~\cite{CLSS}, consisting of freezing the value of the control over a certain duration, called sampling time. The resulting sampled control is then time sparse, by construction.
Therefore, in such a way we obtain a \textit{time sparse and componentwise sparse feedback control}.
\end{remark}

\begin{remark}[Sparse is ``optimal'']
It has been proven in  \cite{CFPT13, CFPT15} that, ``sparse is better" in the following sense:
\begin{center}
\textit{
For every time $t$, $u^\circ(t)$ minimizes $\frac{d}{dt} V(t)$ over all possible feedback controls.
}
\end{center}
In other words, at every instant of time $t$, the above feedback control $u^\circ(t)$ is the best choice in terms of the rate of convergence to flocking.
This means that a policy maker who is not allowed to have prediction on future developments should always consider more favorable to intervene with stronger actions on the fewest possible instantaneous optimal leaders,
rather than
trying to control more agents with minor strength.
\end{remark}

\begin{remark}[Optimal is sparse]
The notion of ``sparsity'' arises naturally in optimal control of multiagent systems. Indeed consider the following simple linear consensus model with an external control
$$
\dot x_{i} = \sum_{j=1^N} a_{ij} (x_{j} - x_{i}) +u_{i}, \qquad i=1,\ldots,N,
$$
where $x_{i} \in \R^{d}$, $a_{ij} \in \R$ and the control $u = (u_{1}, \ldots, u_{N})$ verifies the constraint \eqref{eq:controlbound}
for some given $M >0$. 
Consider the problem of steering the system to the consensus $x_{1}= \dots = x_{N}$ in minimal time $T$. Then the Pontryagin Maximum Principle ensures the existence of a absolutely continuous  nontrivial vector function $p_{1}(t),\ldots,p_{N}(t)$
satisfying the adjoint equation 
$$
\dot p_{i} = \sum_{j=1}^{N} a_{ij} (p_{i} - p_{j})
$$
with final constraint $\sum_{i} p_{i}(T) =0$. The minimality condition for the time optimal control reads
$$
\min \sum_{i=1}^{N} \langle p_{i}(t),u_{i}(t) \rangle,
$$ 
so that the optimal control is, if $\|p_{i}(t)\| > \|p_{j}(t)\|$ for every $j\neq i$,
$$
u_{i}= - M \frac{p_{i}(t)}{\|p_{i}(t)\|} \quad \mbox{ and } u_{j} = 0 \mbox{ for } j\neq i.
$$
In particular the time optimal control is sparse except  when the index for which $\|p_{i}(t)\|$ is maximal is not unique.  However in the generic case in which all interaction coefficients $a_{ij}$ are pairwise distinct then the time optimal control  is sparse for almost every $t \in [0,T]$. 

The previous analysis was done on an approachable toy model. In general, finding the optimal strategy for a consensus model can be very hard. However it is possible to find sparsity features for more complicated optimal control problems, as shown in Section~\ref{sec:optimalcontrolCS} below.
\end{remark}

\subsubsection{Sparse Local Controllability near consensus}

It is possible to show that generically a consensus or an alignment system is controllable near the consensus or alignment manifold by means of a sparse control. More precisely, given a generic pair of configurations sufficiently close to consensus, there exists a strategy, acting on a single agent at every time, that steers the system from one configuration to the other. This property was proven  in \cite{CFPT13,CFPT15} for alignment systems. The proof relies on the fact that given a generic Laplacian matrix $L$ satisfying  some open condition on the coefficients, and  any column vector $B$ with only one component different from $0$, the linear system
$$
\dot x = -L x + B u
$$
verifies the Kalman rank condition for controllability. From this fact it is possible to infer small time local controllability for opinion formation models or alignment models near consensus. Moreover, it is possible to choose the controlled agent a priori. More precisely, the result of \cite{CFPT13,CFPT15} is the following:

\begin{prop}
For almost every consensus there exists a neighborhood in which controllability with time sparse and componentwise sparse control holds.
\end{prop}

This result is easy to establish by linearization around a consensus point. The Kalman condition holds 
for every consensus point verifying and open algebraic condition on the coefficients $a(\|x_{i} - x_{j}\|)$ $i,j=1,\ldots,N$ (whence the ``almost every" of the statement). First order models can be dealt with similarly (see also~\cite{WCB} on opinion formation models).

By stabilization and iterated local controllability along a path of consensus points (note that the set of consensus points is arc-connected), we obtain the following result:

\begin{corollary}
Any point of $(\mathbb{R}^d)^N \times (\mathbb{R}^d)^N$ can be steered to almost any point of the consensus manifold \textit{in finite time} by means of a \textit{time sparse} and \textit{componentwise sparse} control.
\end{corollary}

\subsubsection{Optimal control}\label{sec:optimalcontrolCS}
In~\cite{CFPT13, CFPT15}, the authors have considered, for the finite-dimensional CS model, the optimal control problem with a fixed initial point and free final point, of minimizing the cost functional
$$
\int_0^T \sum_{i=1}^N  \Big( v_i(t) - \frac{1}{N}\sum_{j=1}^N v_j(t) \Big)^2 dt  + \gamma  \sum_{i=1}^N \int_0^T \Vert u_i(t)\Vert  dt
$$
where $\gamma>0$ is fixed, under the constraint $\sum_{i=1}^N  \Vert u_i(t)\Vert \leq M$. As before,
the $\ell^1$-norm (with weight $\gamma$) implies componentwise sparsity features of the optimal control. The proof is done by applying the Pontryagin maximum principle and by developing genericity arguments.

However, because of the coupling between space and velocity, these properties may not be easy to check in practice.
We can note that, if instead of considering the CS model \eqref{eq:CS}, we consider the much simpler HK model \eqref{eq:HK},
then the above optimal control problem becomes quite obvious to analyze. For instance, it is easy to prove that, under generic conditions on the interaction coefficients $a_{ij}$, the optimal control is componentwise sparse.
Such optimal control problems have not yet been considered for the kinetic CS equation \eqref{kineticCS_control}.

Another interesting example of an optimal control problem involves the \textit{collective migration} model \cite{L13, PPS15}, in which the agents (for example migrating birds) aim to align their velocities to a target migration velocity. In this model, not only do the agents interact with each other to evolve as a group as in the CS model, but they also gather clues from the environment to sense the predetermined migration velocity $V$. The control is not an exterior force represented by an additive control as in \eqref{sys_CS_control}. Instead the control is considered to reflect an internal decision making-process between two possible actions: following the group or sensing the target migration velocity. Each agent balances those two forces via a control function $\alpha_i\in [0,1]$. The controlled system writes: 
\begin{equation}\label{Model_Migration}
\begin{cases}
\dot{x_i} = v_i \\
\dot{v_i} = \alpha_i (V-v_i) + (1-\alpha_i)\cfrac{1}{N} \sum\limits_{j=1}^N a(\|x_j-x_i\|) (v_j - v_i)
\end{cases}\quad \text{ for } i\in\{1,...,N\},
\end{equation}
One way to minimize distance from alignment to the target velocity $V$ is to minimize at a given final time the functional $\mathbb{V} = \frac1N \sum_i \|v_i-V\|^2$ with the constraints $0\leq \alpha_i\leq 1$ for all $i \in\{1,...N\}$ and $\sum_i\alpha_i\leq M$. The constraint on the total control strength $\sum_i\alpha_i\leq M$ reflects the fact that it is more energy consuming to sense the target velocity than it is to follow the group, so that only a few individuals can sense it. This naturally divides the group into leaders and followers.
Using the Pontryagin Maximum Principle, Piccoli, Pouradier Duteil and Scharf \cite{PPS15} were able to fully determine an optimal control strategy in the case $M=1$ and $a\equiv 1$. Interestingly, when the initial average velocity $\bar{v}$ is very close to the target velocity $V$, the optimal control strategy requires the existence of an initial \textit{inactivation} time, during which no control should be exerted on the system. This is due to the fact that without control, due to its inherent properties, the system naturally relaxes to alignment of all velocities.

\subsection{Infinite dimension}

%As seen in section \ref{Sec:Overview}, 
%In numerous applications involving risk-taking in economics, pricing models, opinion formation and crowd dynamics, the system is made of a very large number of agents.
%As mentioned above, control of social dynamics systems becomes a very challenging problem when the dimension of the system increases. This is referred to as \textit{the curse of dimensionality}, a term coined by Bellman in the context of dynamic optimization of high-dimensional systems.
%One way around this problem is to move away from the microscopic viewpoint where each agent is considered individually, and consider instead the mean-field limit, which provides a kinetic description of the system (as see in Section \ref{Sec:Overview}). This approach consists of 
%approximating the influence of all individuals on any given one by one averaged effect.
%Numerous ways to control the kinetic equation have been developed.

As in Section \ref{sec_kineticCS}, we want to consider, in some appropriate way, the mean-field limit of the finite-dimensional controlled CS system \eqref{sys_CS_control}. A difficulty comes from the fact that, to ensure existence and uniqueness of the resulting kinetic equation, we need minimal regularity properties of the velocity field, that do not hold true when considering general controls. Another difficulty is that passing to the limit in the previously designed sparse control makes no sense: indeed, in the finite-dimensional model the control has been designed in such a way that, at any instant of time, at most one component of the control is active.
%, meaning that, for instance, one dog is sufficient to control a group of $N$ sheep. 
But when taking the limit $N\rightarrow +\infty$, this is not feasible and the notion of sparsity must be redefined.

\subsubsection{$\Gamma$-limit} \label{Sec:Gamlim}

A first approach in controlling kinetic equations consists of taking the limit of the finite-dimensional controls, in a sense defined by Fornasier and Solombrino \cite{FS14}. 
This combines the concepts of mean-field limit for the probability measure and of $\Gamma$-limit for the control in order to define an appropriate mean-field control for the kinetic equation. 
More specifically, we study the limit when $N\rightarrow \infty$ of the control problem that consists of
%a finite-dimensional control problem in a microscopic setting can be written as 
finding the minimum of the cost functional
\begin{equation}\label{eq:Gamlimmin}
\int_0^T\int_{\R^{2d}}(L(x,v,\mu_N)+\psi(f(t,x,v))d\mu_N(t,x,v)dt
\end{equation}
over all control functions $f$ that satisfy:
\begin{itemize}
\item[(i)] $f:[0,T]\times \R^n\rightarrow\R^d$ is a Carath\'eodory function
\item[(ii)] $f(t,\cdot)\in W_{\text{loc}}^{1,\infty}(\R^n,\R^d)$ for almost every $t\in [0,T]$
\item[(iii)] $|f(t,0)|+\mathrm{Lip}(f(t,\cdot),\R^d)\leq \ell (t)$ for almost every $t\in [0,T]$
\end{itemize} 
where $\ell\in L^q(0,T)$ for a given horizon time $T>0$ and $1\leq q <+\infty$. 
In \eqref{eq:Gamlimmin},
$$
\mu_N(t,x,v)=\frac{1}{N}\sum_{i=1}^N\delta_{(x_i,v_i)}(x,v)
$$
is the atomic measure supported on the phase space trajectories $(x_i(t),v_i(t))\in\R^{2d}$, constrained by satisfying the system: 
\begin{equation}
\begin{cases}
\dot{x}_i = v_i, \\
\dot{v}_i = (H\star\mu_N)(x_i,v_i)+f(t,x_i,v_i)
\end{cases}
\end{equation}
with initial datum $\mu_N^0(t,x,v)=\frac{1}{N}\sum_{i=1}^N\delta_{(x_i^0,v_i^0)}(x,v)$. The notation $H\star\mu_N$ denotes the convolution of $H$ with $\mu_N$, where $H : \R^d \rightarrow \R^d$ is a sublinear and locally Lipschitz continuous interaction
kernel. The functions $L$ and $\psi$ are taken to satisfy appropriate conditions \cite{FS14}.
Applications of this problem include finding ways to influence large crowds, for instance to guide them through an exit in emergency situations. In this context the function $f$ represents an external control on the crowd \cite{FPR14}.

If there exists a compactly supported limit $\mu^0$ to the sequence of atomic measures $\mu_N$ when the number of agents $N$ tends to infinity in the sense of the Wasserstein distance (i.e. $\lim_{N\rightarrow\infty}\mathcal{W}_1(\mu_N^0,\mu^0)=0$), then there exists a subsequence $(f_{N_k})_{k\in\N}$ and a function $f_\infty$ such that $f_{N_k}$ $\Gamma$-converges to $f_\infty$ and $f_\infty$ is a solution of the infinite dimensional optimal control problem 
\begin{equation}
\min_f \int_0^T\int_{\R^{2d}}(L(x,v,\mu_N)+\psi(f(t,x,v))d\mu(t,x,v)dt
\end{equation}
where $\mu$ is the unique weak solution of the kinetic equation
\begin{equation}
\frac{\partial\mu}{\partial t}+\langle v, \nabla_x\mu\rangle = \dv_v \left( (H\star\mu+f)\mu \right)
\end{equation}
with initial datum $\mu^0$.

\subsubsection{Control by leaders} \label{Sec:ContLead}

A common way to control a large crowd is to act on a selected few individuals that will behave as leaders to guide the crowd. 
In the case of finite dimensional systems, it is frequent to look for controls that are vanishing for most of the agents and for most of the time. 
These strategies are referred to as sparse control, as seen in Section \ref{Sec:ControlFinite}. Such controls have obvious advantages, being both moderate in the external control and parsimonious in the number of agents controlled. 
Extending the concept of leaders to mean-field limits is not straightforward. Indeed, when representing the crowd with a particle density distribution, the action of a finite number of agents becomes negligible compared to the size of the crowd. Albi and Pareschi have looked at the microscopic-macroscopic limit of such systems, see \cite{AP13}. In \cite{FPR14}, Fornasier, Piccoli and Rossi solve this problem by using a mixed granular-diffuse description of the crowd and prove convergence of the solution of the finite-dimensional problem when the number of followers tends to infinity to the solution of this new system. Similar approaches involving the coupling of microscopic dynamics for the leaders and macroscopic dynamics for the followers were adopted by  Albi, Bongini, Cristiani and Kalise in \cite{ABCK16} and by Colombo and Pogodaev in \cite{CP13}. Here we present the approach of \cite{FPR14}. Let $(y_k,w_k)$ denote the space-velocity variables of the $m$ leaders of the crowd, and $(x_i,v_i)$ those of the $N$ followers, so that for a given locally Lipschitz interaction kernel with sublinear growth $H$, 
\begin{equation}\label{eq:leadfollfin}
\begin{cases}
\dot{y}_k = w_k \\
\dot{w}_k = H \star \mu_N(y_k,w_k)+H\star\mu_m(y_k,w_k) + u_k \quad \text{for } k\in\{1,...,m\} \\
\dot{x}_i = v_i \\
\dot{v}_i = H \star \mu_N(x_i,v_i)+H\star\mu_m(x_i,v_i) \quad \text{for } i\in\{1,...,N\} 
\end{cases}
\end{equation}
where
\begin{equation}
\mu_N(t,x,v)=\frac{1}{N}\sum_{i=1}^N\delta_{(x_i(t),v_i(t))} \text{ and }
\mu_m(t,y,w)=\frac{1}{m}\sum_{k=1}^m\delta_{(y_k(t),w_k(t))}
\end{equation}
and $u_k:[0,T]\rightarrow\R^d$ are measurable controls for $k\in\{1,...,m\}$.
The optimal control problem consists of finding solutions of 
\begin{equation}\label{eq:contfin}
\min_{u_k} \int_0^T(L(y(t),w(t),\mu_N(t))+ \frac{1}{m}\sum_{k=1}^m |u_k(t)| ) dt.
\end{equation} 
where $(y(t),w(t),\mu_N(t))$ are subject to the dynamics \eqref{eq:leadfollfin}.

The authors of \cite{FPR14} showed that a mean-field limit of system \eqref{eq:leadfollfin} when $N$ tends to infinity can be derived as the coupling of controlled ODE's for the evolution of the leaders' positions and velocities and of a PDE for the compactly supported probability measure $\mu$ of the followers in the position-velocity space: 
\begin{equation}\label{eq:leadfollinf}
\begin{cases}
\dot{y}_k = w_k \\
\dot{w}_k = H \star (\mu+\mu_m)(y_k,w_k) + u_k \quad \text{for } k\in\{1,...,m\} \\
\partial_t\mu + \langle v, \nabla_x\mu\rangle = \dv_v ((H\star (\mu+\mu_m))\mu).
\end{cases}
\end{equation}

Moreover, the optimal controls $u_N^*$ of the finite dimensional control problem \eqref{eq:contfin}-\eqref{eq:leadfollfin} converge weakly for $N\rightarrow\infty$ to optimal controls $u^*$ that are solutions of: 
\begin{equation}\label{eq:continf}
\min_{u_k} \int_0^T(L(y(t),w(t),\mu(t))+ \frac{1}{m}\sum_{k=1}^m |u_k(t)| ) dt.
\end{equation}
where $(y(t),w(t),\mu(t))$ are subject to the dynamics \eqref{eq:leadfollinf}.
In \cite{FPPR14}, these results were extended to Mayer-type minimization problems. 
Such formulations are particularly useful in that they combat the \textit{curse of dimensionality}. Indeed, even though the total number of agents is allowed to tend to infinity, the number of controlled ones stays bounded and small, which keeps the numerical computations feasible. 

\subsubsection{Controlled kinetic Cucker-Smale model}

The first two approaches to controlling kinetic systems reported in Sections \ref{Sec:Gamlim} and \ref{Sec:ContLead} consist of finding an optimal control for the finite-dimensional system and passing to the limit (in some appropriate sense) when the number of agents tends to infinity.
Another approach involves controlling the PDE directly. This was done for instance by Piccoli, Rossi and Tr\'elat \cite{PRT15}.
%As in Section \ref{sec_kineticCS}, we want to consider, in some appropriate way, the mean-field limit of the finite-dimensional controlled Cucker-Smale system \eqref{sys_CS_control}. A difficulty is coming from the fact that, to ensure existence and uniqueness of the resulting kinetic equation, we need minimal regularity properties of the velocity field, that do not hold true when considering general controls. Another difficulty is that passing to the limit the previously designed sparse control makes no sense: indeed, in the finite-dimensional model the control has been designed in such a way that, at any instant of time, at most one component of the control is active, meaning that, for instance, one dog is sufficient to control a group of $N$ sheep. But when taking the limit $N\rightarrow +\infty$, this will make no sense.
A difficulty arises when the control is componentwise sparse (see Section \ref{Sec:Control-Finite-Sparse}). Keeping just one component of the control active is only practical in finite dimension. One way to translate this criterion to infinite-dimensional problems is to 
 pass to the limit by keeping proportions: given some fixed $c\in(0,1)$,  assume that the control acts on $cN$ agents of the $N$-sized group.
 %, in order to control a group of $N$ sheep, we have (around) $cN$ dogs. 
 It is then easy to see how to modify the sparse control designed for the finite-dimensional model in order to fit this new requirement, and it makes sense to pass to the limit $N\rightarrow +\infty$. The real number $c$ represents the proportion of the crowd on which one is allowed to act.
When passing to the limit, we obtain a control domain, denoted $\omega(t)$ in what follows, representing the controlled part of the crowd at time $t$.

Let us now be more precise.
Following \cite{PRT15}, we consider the mean-field limit
\begin{equation}\label{kineticCS_control}
\partial_t \mu +\langle v, \nabla_x \mu\rangle+\mathrm{div}_v \left( (\xi[\mu] +\chi_\omega u) \mu \right)=0 
\end{equation}
with $u\in L^\infty(\R\times\R^d\times\R^d)$ and $\omega(t)\subset\R^d$ measurable, such that
\begin{equation}\label{cont_u}
\Vert u(t,\cdot,\cdot)\Vert _{L^\infty(\R^d\times\R^d)}\leq 1,
\end{equation}
standing for a bounded external action, and
\begin{equation}\label{cont_mu}
\mu(t)(\omega(t)) = \int_{\omega(t)} d\mu(t)(x,v) \leq c,
\end{equation}
modeling that one is allowed to act only on \textit{a given proportion $c$ of the crowd}.
This is the natural notion of \textit{sparse control} in the infinite-dimensional setting.

A possible variant is to consider the constraint
\begin{equation}\label{cont_om}
\vert\omega(t)\vert=\int_{\omega(t)} \, dx\, dv\leq c ,
\end{equation}
representing a limit on the space of configurations.

Hence, now the control is $\chi_\omega u$, and consists of choosing, at any instant of time, the control domain $\omega(t)$, and the force $u(t,\cdot,\cdot)$ with which one acts on the crowd along the control domain.
Note that it is not so common in the existing literature to control not only an external force but also a domain.

\begin{theorem} \cite{PRT15}
For every $\mu^0\in\mathcal{P}_c^{ac}(\R^d\times\R^d)$, there exists a control $\chi_\omega u$ satisfying \eqref{cont_u} and either \eqref{cont_mu} or \eqref{cont_om}, and the corresponding solution $\mu\in C^0(\R; \mathcal{P}_c^{ac}(\R^d\times\R^d))$ such that $\mu(0)=\mu^0$ converges to consensus.
\end{theorem}

\begin{remark}
The strategy to prove this theorem is quite long and technical, and is not reported in detail here. We just give hereafter the main intuitive idea. Writing the controlled kinetic CS equation \eqref{kineticCS_control} as
\begin{equation}\label{CS_velocityfield}
\partial_t \mu + \mathrm{div}_{(x,v)} \left( V_{\omega,u}[\mu]  \mu \right) = 0 ,
\end{equation}
with the controlled velocity field
$$
V_{\omega,u}[\mu] = \begin{pmatrix} v \\\xi[\mu]+ \chi_\omega u \end{pmatrix} ,
$$
the \textit{controlled particle flow} $\Phi_{\omega,u}(t)$ generated by $V_{\omega,u}[\mu(t)]$ yields the characteristics
$$
\dot x(t) = v(t),\qquad \dot v(t) = \xi[\mu(t)](x(t),v(t)) + \chi_{\omega(t)} u(t,x(t),v(t)).
$$
This is a control system, describing the (controlled) motion of particles. As in the uncontrolled case, the measure, solution of the kinetic equation, is then the pushforward of the initial measure:
$$
\mu(t)=\Phi_{\omega,u}(t)\#\mu^0 .
$$
Having these facts in mind, we adopt, as in the finite-dimensional case, a \textit{shepherd control design strategy}: at every instant of time, we choose $\omega(t)$ and $u(t)$ such that the controlled velocity field $V_{\omega,u}[\mu(t)]$ points inwards the invariant domain, thus confining the population. This implies that the size of $\mathrm{supp}_v(\mu(t))$ (velocity support of the measure) decreases exponentially in time.
This construction is carried out in \cite{PRT15}, piecewise in time, and in an algorithmic way, thus resulting in an explicit control strategy such that
\begin{itemize}[nolistsep]
\item$\omega(t)$ is piecewise constant in $t$,
\item $u(t,x,v)$ is piecewise constant in $t$ for $(x,v)$ fixed, $C^0$ and piecewise linear in $(x,v)$ for $t$ fixed.
\end{itemize}
An important difficulty in dealing with the kinetic equation \eqref{CS_velocityfield} is to ensure existence and uniqueness of the solution. Indeed, for general controls the velocity field is not regular enough. An essential feature of the control strategy designed in \cite{PRT15} is that the control is piecewise smooth, and then existence and uniqueness of the solution are ensured in an iterative way. Actually, the solution $\mu$ of \eqref{CS_velocityfield} remains absolutely continuous and of compact support, and it becomes singular only in infinite time.

Note that, as in the finite-dimensional setting, the control is switched-off when entering the consensus region: given any $\mu^0\in\mathcal{P}_c^{ac}(\R^d\times\R^d)$, there exists $T(\mu^0)\geq 0$ such that $u(t,x,v)=0$ for every $t>T(\mu^0)$. Then the solution reaches consensus (here, a Dirac mass) in infinite time, and remains absolutely continuous in-between.
\end{remark}

Another variant is not to impose a constraint on the control but to penalize its spread in the cost function, as done in \cite{FS14}. More specifically, for a PDE constrained problem
\begin{equation}
\partial_t\mu+\langle v,\nabla_x\mu\rangle = \dv_v((H\star\mu+f)\mu),
\end{equation}
the control $f$ is the minimizer of the chosen cost: 
\begin{equation}
\mathcal{E}_\psi(f):=\int_0^T\int_{\R^2d} (L(x,v,\mu) +\psi (f(t,x,v)))d\mu(t,x,v)dt
\end{equation}
where a relevant choice for $\psi$ is for instance $\psi(\cdot):=\gamma |\cdot |$ for $\gamma>0$. This promotes the sparsity of $f$ thanks to the $\ell^1$ norm penalization.

\subsubsection{Boltzmann-type control for consensus dynamics}

As stressed in the previous sections, controlling a large number of agents is computationally expensive, and even sometimes unfeasible. For instance, the Pontryagin Maximum Principle applied to the minimization problem 
\begin{equation}
\min_{u(t)\in[u_L,u_R]} \int_0^T \frac{1}{N}\sum_{j=1}^N(\frac{1}{2}(x_j-x_d)^2+\frac{\kappa}{2}u^2)ds,
\end{equation} 
where $x_d$ is the desired state, for the controlled system
\begin{equation}
\dot{x}_i=\frac{1}{N}\sum_{j=1}^N a(x_i,x_j)(x_j-x_i)+u
\end{equation}
requires solving the equation for the adjoint vector backwards in time over the whole interval $[0,T]$, which is extremely costly for a large number of agents. In \cite{AHP15}, Albi, Herty and Pareschi develop an alternative approach consisting of solving the control problem on a sequence of reduced time-intervals. This iterative method is called \textit{model predictive control}.
The horizon-receding strategy allows to embed the minimization of the cost functional into the particle interactions. 

As done in \cite{AHP15}, let $I = [-1,1]$ represent a bounded set of opinions such that $x_i(t)\in I$, $i=1,...,N$. Alternatively, in the multidimensional setting, one can take $I=\Sp^d$.
Consider the case of binary Boltzmann dynamics with two interacting agents $i$ and $j$, then their positions $x^{n+1}$ at time $(n+1)\Delta t$ depends on the previous state in the following way: 
\begin{equation}
\begin{cases}
x_i^{n+1} = x_i^n + \frac{\Delta t}{2} a(x_i^n,x_j^n)(x_j^n-x_i^n)+\frac{\Delta t}{2} U(x_i^n,x_j^n) \\
x_j^{n+1} = x_j^n + \frac{\Delta t}{2} a(x_i^n,x_j^n)(x_j^n-x_i^n)+\frac{\Delta t}{2} U(x_j^n,x_i^n)
\end{cases}.
\end{equation} 
The model predictive control performed on a single prediction time-horizon allows the explicit expression: 
%$\Delta t \; u^n = \frac{\Delta t}{2(\nu+\Delta t)}U_\alpha$, where $\alpha = \Delta t/2$ and $U_\alpha$ is given by: 
\begin{equation}
\eta U(x_i^n,x_j^n) = \frac{\beta}{2} \left( (x_d-x_j^n)+(x_d-x_i^n) + \eta (a(x_i,x_j)-a(x_j,x_i))(x_j^n-x_i^n)\right)
\end{equation}
where $\beta:=\frac{2\eta}{\kappa+2\eta}$ and $\eta = \Delta t/2$.
The kinetic Boltzmann equation is obtained by introducing the density distribution of particles $\mu(x,t)$ belonging to the space of probability measures. Then two agents $x$ and $y$ modify their states according to: 
\begin{equation}\label{eq:binint}
\begin{cases}
x^* = x + \eta \left( a(x,y)(y-x) +  U(x,y) \right) \\
y^* = y + \eta \left( a(y,x)(x-y) +  U(y,x) \right) \\
\end{cases}
\end{equation} 
%where $\alpha$ and $\beta:=\frac{2\eta}{\kappa+2\eta}$ are the strength of the compromise and of the control.
% Notice that two additional noise terms have been added, characterized by the random variables $\Theta_1$ and $\Theta_2$, and by the function $0\leq D(w)\leq 1$ that represents the local relevance of the opinion spreading due to diffusion.

For a test function $\phi(x)$, we write: 
\begin{equation}\label{eq:KinBoltzmann}
\frac{d}{dt}\int_I\phi(x)\mu(x,t)dx = \lambda \int_{I^2}\left(\phi(x^*)+\phi(y*)-\phi(x)-\phi(y)\right)\mu(x,t)\mu(y,t)dxdy
\end{equation}
where $\lambda$ represents a constant rate of interaction and we considered that $a(x,y)=a(y,x)$. 
%here $\langle\cdot\rangle$ represents the expectation with respect to $\Theta_1$ and $\Theta_2$, and
% $B_\mathrm{int}$ is a nonnegative interaction kernel containing the information from the binary interaction \eqref{eq:binint} such that $B_\mathrm{int}(x,y) = \eta \chi(|x^*|\leq 1) \chi(|y^*|\leq 1)$ (with $\chi(\cdot)$ the indicator function).
This allows us to show that the limit of the average position $m_\infty:=\lim_{t\rightarrow\infty}m(t)$ where $m(t) = \int_Ix\mu(x,t)dx$ stays close to the desired state $x_d$, and in the symmetric case $a(x,y)=a(x,y)$, we even have $m_\infty=x_d$. Moreover, if the interaction kernel is simplified to $a(x,y)=1$, then one shows that the particle distribution converges to the Dirac measure $\delta(x-x_d)$ centered in the desired state $x_d$, which implies that the system reaches consensus. 

To derive the asymptotic limit of the model while retaining the memory of the binary interactions \eqref{eq:binint}, one can refer to the so-called 
\textit{quasi-invariant opinion limit} \cite{T06}. This consists of adapting to the context of consensus models the concept of
\textit{grazing collision limit} used to consider long time solutions of the Boltzmann equation \cite{V98}. For a summary of these concepts, see \cite{PT13}. 
Here, this is done by rescaling time in \eqref{eq:KinBoltzmann}.
%, which consists of scaling parameters appropriately. 
In particular, taking $\eta=\epsilon$, $\lambda = 1/\epsilon$,
% $\nu = \epsilon\kappa$, 
the limit when $\epsilon$ tends to zero leads to a controlled kinetic equation of type
$$
\mu_t = \dv_x ( (\xi[\mu]+\zeta[\mu])\mu), 
$$  
where 
$$
\xi[\mu](x) = \int_I a(x,y)(y-x)\mu(y)dy
$$
and
$$
\zeta[\mu](x) = \int_{\R^d} K(x,y) d\mu(y), \text{ for } K(x,y) =\frac{1}{\kappa}((x_d-x)+(x_d-y)).
$$
This approach can be easily extended to the case of the CS model and leader-follower model, see \cite{AP16,APZ14}.

\subsubsection{Mean-field games}

Another common way to deal with control of large systems is to use \textit{mean-field games}, a theory that was introduced by Lasry and Lions in 2006 \cite{GLL11,LL07} and independently by Caines \cite{C13,HMC06}.
A wealth of results have since then been obtained 
in this game-theoretic setting,  
 considering that each agent makes a decision in order to optimize a given cost based on its available information, (see Degond, Liu and Ringhofer \cite{DLR13}). 
 For instance, in applications to economics, 
it is meaningful to study Nash equilibria, a stable state in which no agent can improve its cost by changing alone its strategy.

Mean-field games are used to consider each agent's individual decision, given his knowledge of the system. Consequently, most applications of mean-field games are found in economics, where each agent strives to optimize its wealth given the current state of the market (Gu\'eant, Lasry and Lions \cite{GLL11}). For instance, price formations models can be derived by dividing the population into buyers and sellers. Other applications involve crowd modeling (Lachapelle and Wolfram \cite{LW11}).

In \cite{BP13, BP14}, Bardi and Priuli derive the mean-field limit of a stochastic differential game with $N$ players: 
\begin{equation}\label{eq:stoch}
dX_t^i = (A^iX_t^i - \alpha_t^i)dt + \sigma^idW_t^i, \quad X_0^i\in \R^d, \quad i=1,...,N.
\end{equation}
where $A^i$ and $\sigma^i$ are $d\times d$ matrices, $(W_t^i)_{i\in\{1,...,N\}}$ are $N$ independent $d$-dimensional standard Brownian motions, and $\alpha_t^i$ is a process adapted to $W_t^i$ representing the control of the $i$-th player, designed to minimize the quadratic running cost
\begin{equation}
J^i(X,\alpha^1,...,\alpha^N):=\liminf_{T\rightarrow\infty}\frac{1}{T}\mathbb{E}\left[ \int_0^T\frac{(\alpha_t^i)^T R^i \alpha_t^i}{2} + (X_t-\bar{X}_i)^T Q^i (X_t-\bar{X}_i)dt\right].
\end{equation}
Here $\mathbb{E}$ denotes the expected value, $R^i$ are positive definite $d\times d$ matrices, $Q^i$ are symmetric $Nd\times Nd$ matrices and $\bar{X}_i$ are reference positions. Following the approach of \cite{LL07}, the corresponding system of $N$ nonlinear Hamilton-Jacobi-Bellman PDEs coupled with $N$ Kolmogorov-Fokker-Plank equations can be derived:
\begin{equation}\label{eq:HJB-KFP}
\begin{cases}
-\tr (\nu^iD^2v^i)+\mathcal{H}^i(x,\nabla v^i)+\lambda^i = f^i(x;m^1,...,m^N) \\
-\tr((\nu^iD^2 m^i)-\dv (m^i\frac{\partial \mathcal{H}^i}{\partial p}(x,\nabla v^i))=0 \\
\int_{\R^d}m^i(x)dx = 1, \quad m^i>0,
\end{cases}
\end{equation}
where the unknowns are the functions $v^i$, the scalars $\lambda^i$ and the measures $m^i$, and $\mathcal{H}^i$ denotes the $i-$th Hamiltonian.
This leads to Nash equilibria obtained by affine feedback.
% for the invariant measure of the process associated with the Nash equilibrium.
In order to derive the mean-field limit of \eqref{eq:stoch}, we need to consider that the players are \textit{nearly identical}, that is that they are influenced in the same way by pairs of other players, have the same control systems ($A^i=A$), the same costs of the controls ($R^i=R$), the same reference positions ($\bar{X}_i=X_d$) and the same primary costs of displacement. Then given suitable conditions, there exist unique solutions to the system of HJB-KFP equations \eqref{eq:HJB-KFP}. 
Moreover, those solutions converge when $N\rightarrow\infty$ to the solutions of a mean-field system of the form:
\begin{equation}\label{eq:HJB-KFP_inf}
\begin{cases}
-\tr (\nu D^2v)+\nabla v^T \frac{R^{-1}}{2}\nabla v - \nabla v^T Ax+\lambda = \hat{V}[m](x) \\
-\tr((\nu D^2 m)-\dv (m\cdot R^{-1}\nabla v -Ax) =0 \\
\int_{\R^d}m(x)dx = 1, \quad m>0.
\end{cases}
\end{equation}

Mean-field games consist of optimizing control strategies over a large time horizon. This is both computationally expensive and not always realistic, since individual agents might not have access to information on the state of the system in the distant future. Instead, some strategies called \textit{best reply strategies} compute the best instantaneous response given the present state of the system, by steepest gradient descent \cite{DLR13}. This type of control is suboptimal in the long term, but is in some ways more realistic and more feasible. A good compromise between standard large-time mean-field approaches and best reply strategies consists of minimizing the cost function over a small shifting time horizon \cite{AHP15}, as done in \textit{Model Predictive Control}, see Maciejowski, Goulart and Kerrigan \cite{M07} and Degond, Herty and Liu \cite{DHL16}.

\section{Generalizations}

There exist many possible generalizations of the models considered above. One way to better adapt the model to the phenomenon of interest is to use general interaction potentials. For example, in the case of animal group modeling, Carillo et al propose to take into account the cone of vision of each animal $i$ to define its influential set of neighbors $\NN_i$ \cite{CFTV10}. This naturally singles out a certain number of instantaneous leaders, defined as the animals whose cones of vision are pointed outwards so that they do not follow any other agent. Such dynamics are expected to lead to clustering of the group into a finite number of subgroups each following a leader.
 
In \cite{CFP11}, Cristiani, Frasca and Piccoli studied the effect of anisotropic interaction regions on the shape of the group. Around each agent are defined are zone of attraction and a zone of repulsion, that can each be isotropic or anisotropic. Depending on the nature of those interactions, simulations show that various patterns can be obtained in the group: crystal-like clusters of individuals, lines or V-like formations.

In \cite{DCBC06}, D'Orsogna, Chuang, Bertozzi and Chayes propose a model to take into account self-propelling, friction and attraction-repulsion effects. More specifically, each agent's velocity is defined as $\dot v_i = (\alpha -\beta\|v_i\|^2)v_i - \nabla_i U(x_i)$, where $U(x_i)$ is 
 the Morse potential, used to include attractive and repulsive ranges. For agent $i$, $U(x_i)=\sum_{j\neq i} (C_r e^{-\|x_i-x_j\|/l_r} + C_a e^{-\|x_i-x_j\|/l_a})$ where $l_r$ and $l_a$ are respectively the repulsive and attractive ranges and $C_r$ and $C_a$ the repulsive and attractive amplitudes. In the case of animal groups or other biological applications, the most relevant cases are $l_a>l_r$ and $C_a>C_r$, for short-range repulsion and long-range attraction.
 The parameter $\alpha$ models the self-propulsion capacity of agent $i$, while the parameter $\beta$ represents a friction according Rayleigh's law. This model gives rise to various types of regimes, depending on the ratios $l_r/l_a$ and $C_r/C/a$. The system is said to be \textit{H}-stable if the total potential energy is bounded below by a multiple of the number of agents $N$, i.e. $U\geq -BN$ for some constant $B\geq 0$. \textit{H}-stability ensures that the system does not collapse when $N\rightarrow\infty$, so that the particles form a crystal-like structure. When the system is not \textit{H}-stable, it is said to be \textit{catastrophic}, and as $N\rightarrow\infty$, the inter-particle intervals shrink to zero. 
 
Another common generalization of the models described in this chapter consists of adding  white noise to the dynamics. For example, in \cite{YEECBKIMS09}, Yates et al. argue that locusts use white noise to maintain swarm alignment. This claim is supported by experimental evidence, itself modeled using the kinetic Fokker-Planck equation with noise. See also \cite{HLL09} for stochastic models.

%Many open control problems range over cluster control, control of opinion formation, black swan problems, cancerology, robot flotilla formation, etc. 

%There are many possible generalizations of the models considered above, among which:
%\begin{itemize}
%\item general potentials (friction, attraction/repulsion, ...), cone vision constraints, leadership (\cite{WCB,WCKB}), ...: 
%Cucker, Dong, Ha, Ha, Kim, Leonard, Motsch, Slemrod, Tadmor.
%\item stochastic aspects (adding noise): Carrillo, Cucker, Fornasier, Ha, Lee, Levy, Mordecki, Toscani.
%\item other models in infinite dimension (hydrodynamic, kinetic, mean-field limit):
%Carrillo, Degond, Fornasier, Ha, Hascovek, Motsch, Rosado, Tadmor, Toscani.
%\end{itemize}